\documentclass[12pt]{article}

\RequirePackage[OT1]{fontenc}
\RequirePackage[amsthm,amsmath,natbib,noinfoline]{imsart}
\RequirePackage[dvips]{hyperref}
\RequirePackage{hypernat}

\startlocaldefs
\numberwithin{equation}{section}
\theoremstyle{plain}

\endlocaldefs

\usepackage[mathscr]{eucal}
\usepackage{graphicx}
\usepackage{latexsym}
\usepackage{amssymb}
\usepackage{psfrag}
\usepackage{epsfig}

\newtheorem{theorem}{Theorem}[section]
\newtheorem{remark}[theorem]{Remark}
\newtheorem{ex}[theorem]{Example}
\newtheorem{lemma}[theorem]{Lemma}
\newtheorem{cor}[theorem]{Corollary}
\newtheorem{defn}[theorem]{Definition}

\newtheorem{prop}[theorem]{Proposition}

\newcommand{\sm}{\setminus}

\newcommand{\noi}{\noindent }
\newcommand{\ba}{\begin{array}}
\newcommand{\ea}{\end{array}}
\newcommand{\bea}{\begin{eqnarray}}
\newcommand{\eea}{\end{eqnarray}}
\newcommand{\beas}{\begin{eqnarray*}}
\newcommand{\eeas}{\end{eqnarray*}}
\newcommand{\be}{\begin{equation}}
\newcommand{\ee}{\end{equation}}
\newcommand{\bt}{\begin{theorem}}
\newcommand{\et}{\end{theorem}}
\newcommand{\bc}{\begin{center}}
\newcommand{\ec}{\end{center}}
\newcommand{\ben}{\begin{enumerate}}
\newcommand{\een}{\end{enumerate}}

\newcommand{\beq}{\begin{eqnarray*}}
\newcommand{\eeq}{\end{eqnarray*}}
\newcommand{\ei}{\end{itemize}}

\newcommand{\bes}{\begin{equation*} }
\newcommand{\ees}{\end{equation*} }

\newcommand{\ds}{\displaystyle}

\newcommand{\ve}{\varepsilon}

\def\Z{\mathbb{Z}}
\def\R{\mathbb{R}}

\def\N{\mathbb{N}}
\def\ind{\mathbb{I}}

\newcommand{\ra}{\rightarrow}

\newcommand{\cspace}{{\cal C}\left[ 0,\infty\right)}

\newcommand{\dspace}{{\cal D}\left[ 0, \infty\right)}
\newcommand{\dspacepl}{{\cal D}^+\left[ 0, \infty\right)}
\newcommand{\dspacemin}{{\cal D}^-\left[ 0, \infty\right)}

\def\gzero{\Gamma_0}

\def\bgamma{\overline{\Gamma}}
\def\bGamma{\overline{\Gamma}}

\def\bphi{\bar{\phi}}

\def\MF{{\mathcal F}}

\def\MV{{\mathcal V}}

\def\P{\mathbb{P}}
\def\E{\mathbb{E}}

\def\diff{\Delta}

\def\ar{r}
\newcommand{\gt}{[\ell (\cdot), \ar (\cdot)]}
\def\ltime{Y}

\def\eps{\varepsilon}

\input{epsf.tex}

\begin{document}

\begin{frontmatter}
\title{The Skorokhod problem in a time-dependent interval}
\runtitle{The Skorokhod problem in a time-dependent interval}

\begin{aug}
\author{\fnms{Krzysztof} \snm{Burdzy}\thanksref{t1}\ead[label=e1]{burdzy@math.washington.edu}},
\author{\fnms{Weining} \snm{Kang}\ead[label=e2]{weikang@andrew.cmu.edu}}
\and
\author{\fnms{Kavita}
  \snm{Ramanan}\thanksref{t3}\ead[label=e3]{kramanan@math.cmu.edu}}

\thankstext{t1}{Partially supported by NSF Grant DMS-0600206.}
\thankstext{t3}{Partially supported by NSF Grants DMS-040691, DMS-0405343.}
\affiliation{University of Washington and Carnegie Mellon University}
\date{today}
\end{aug}

\begin{abstract}
We consider the Skorokhod problem in a time-varying interval.
We prove existence and uniqueness for the solution.  We also
express the solution in terms of  an explicit formula.
Moving boundaries may generate singularities when they touch.
We establish two sets of sufficient conditions on the moving boundaries
that guarantee that the variation of the local time of the associated
reflected Brownian motion is, respectively, finite and infinite.
We also apply these results to study the semimartingale property of
a class of two-dimensional reflected Brownian motions.
\end{abstract}



\end{frontmatter}

\section{Introduction}

We consider the Skorokhod problem with two moving boundaries.
Informally speaking, the problem is concerned with reflecting or
constraining a given path in a space-time region defined by two moving
boundaries. We will address several problems inspired by recent
related developments.  First, we  study the question of existence
and uniqueness to a slight generalization of the Skorokhod problem, which
we refer to as the extended Skorokhod problem.
We  show that the
solution not only exists and is unique, but can be represented in
terms of an  explicit and rather simple  formula.
Second, we  prove some monotonicity relations for solutions to the
extended Skorokhod problem.
Similar monotonicity properties are quite obvious when there is only
one reflecting boundary; they are not so obvious in our context.
Moreover, we  study the issue of whether the local time of the reflected Brownian
motion has finite or infinite total variation. This issue  arises
when the boundary of the domains are allowed to meet and
 is related to the question of whether the reflected Brownian motion
is a semimartingale (see, for example, \cite{kanram1,wil}).
Finally, we  apply our analysis of one-dimensional reflected
Brownian motion in a time-dependent domain to study the
behavior of the local time and, in particular, the
semimartingale property of a class of two-dimensional reflected
Brownian motions in a fixed domain that were studied in
\cite{burtoby,ram1,ramrei2,wil}. Reflecting Brownian motions in
time-dependent domains arise in queueing theory \cite{manmas},
statistical physics \cite{burnua,souwer}, control theory
\cite{elkar-karatzas} and finance \cite{elkar_kapo}.

The present paper is related to several articles. First, the
papers \cite{krulehramshr2} and \cite{krulehramshr1} present an
explicit formula for the Skorokhod mapping in the simpler
setting of a constant interval $[0,a]$. Second, the works
\cite{burchesyl04} and \cite{burchesyl03} contain an analysis
of Brownian motion reflected on one moving boundary. In
particular, the second paper presents results on singularities
at rough boundary points. In the present paper, we analyze
singularities due to the interaction of two moving boundaries.
In our context, a ``singularity'' means the infinite variation
of the local time process. Finally, the paper  \cite{burtoby}
(see also \cite{wil}) studied a special case of two-dimensional
reflected Brownian motion in thorn-like domains, with all
reflection vectors parallel to the same straight line. We
establish the somewhat suprising result that this
two-dimensional reflected Brownian motion is not a
semimartingale, irrespective of the particular shape of the
thorn.  In addition, we also provide new proofs of some of the
qualitative results of the papers \cite{burtoby,wil}.

The rest of the paper is organized as follows. We start with a
short section collecting the notation used throughout the
paper. Section \ref{sec-sptv} is devoted to the foundational
results---existence, uniqueness and an explicit formula for the
so-called extended Skorokhod mapping. Section \ref{sec-comp}
contains some ``comparison'' or ``monotonicity'' results.
Finally, Sections \ref{subs-lower} and \ref{subs-upper} present
theorems on the local time of reflected Brownian motion in a
time-dependent interval.  These results are applied in Section
\ref{subs-2drbm} to study the local time of a class of
two-dimensional reflected Brownian motions.

\subsection{Notation}
 \label{subs-notat}

We use $\dspace$ to denote the space of c\`{a}dl\`{a}g
functions (i.e., continuous on the right with finite left
limits) that are defined on $[0,\infty)$ and take values in
$(-\infty,\infty)$. The space of c\`{a}dl\`{a}g functions
taking values in $[-\infty, \infty)$ (respectively,  $(-\infty,
\infty]$) will be denoted $\dspacemin$ (respectively,
$\dspacepl$). Given two functions $f \in \dspacemin$, $g \in
\dspacepl$, we will say $f \leq g$ (respectively, $f < g$) if
$f(t) \leq g(t)$ (respectively, $f(t) < g(t)$) for every $t \in
[0,\infty)$. We let $\cspace$ represent the subspace of
continuous functions in $\dspace$. We denote the variation of a
function $f$ on $[t_1,t_2]$ by $\MV_{[t_1,t_2]}(f)$. We denote
by $\ell (\cdot)$ a generic function in $\dspacemin$ and by
$\ar (\cdot)$ a generic function in $\dspacepl$, and assume
that $\ell \leq \ar$.

Moreover, given $a, b \in \R$, denote $a\wedge
b\doteq\min\{a,b\}$, $a\vee b\doteq\max\{a,b\}$, and $a^+\doteq
a\vee 0$.  We denote by $\ind_A$ the indicator function of a
set $A$.

We also use the following abbreviations, whose meaning will be
explained later: SP---Skorokhod problem, SM---Skorokhod map,
ESP---extended Skorokhod problem, ESM---extended Skorokhod map,
BM---Brownian motion, RBM---reflected Brownian motion.

\section{Skorokhod and Extended Skorokhod Maps in a Time-Dependent Interval}
\label{sec-sptv}

The so-called Skorokhod Problem (SP) was introduced in
\cite{sko2} as  a convenient tool for the construction of
reflected Brownian motion (RBM) in the {\em time-independent}
domain $[0,\infty)$. Specifically, given a function $\psi \in
\dspace$, the SP on $[0,\infty)$ consists of identifying a
non-negative function $\phi$ such that the function $\eta
\doteq \phi - \psi$ is non-decreasing and, roughly speaking, increases only at
times $t$ when $\phi (t) = 0$. It was shown in \cite{sko2} that
there is a unique mapping that takes  any given $\psi \in
\cspace$ to the corresponding function $\phi$ (the extension to
$\psi \in \dspace$ is straightforward). Moreover, this mapping,
which we shall refer to as the Skorokhod map (SM) on
$[0,\infty)$ and denote by $\gzero$, admits the explicit
representation
 \be \label{def-g0}
 \gzero (\psi )(t) = \psi (t)
 + \sup_{s \in [0,t]} \left[ -\psi(s)
 \right]^+,\quad\psi\in\dspace.
 \ee
Given  a Brownian motion  (BM) $B$ on $\R$ with $B(0) = 0$, and
any $x \geq 0$,
 the process $W = \gzero (x + B)$ defines RBM on $[0,\infty)$,
starting at $x$. More generally, due to the Lipschitz
continuity of the map $\Gamma_0$, standard Picard iteration
techniques can be used to construct solutions to stochastic
differential equations with reflection on $[0,\infty)$, under
the usual Lipschitz assumptions on the drift and diffusion
coefficients.

In a similar fashion, the generalizations of the SP given in Section \ref{subs-spdef} will be
the basis for the construction of $1$-dimensional RBM in a time-dependent interval.
We also establish some basic properties of these generalizations in
Section \ref{subs-spdef} and then provide an explicit formula for the ESM in Section \ref{subs-explicit}.

\subsection{Basic Definitions and Properties}
\label{subs-spdef}

We first describe the SP on a time-varying interval $\gt$.

\begin{defn} {\bf (Skorokhod problem on $[\ell(\cdot),r(\cdot)]$)}
\label{def-smap1} Suppose that $\ell  \in \dspacemin$, $\ar \in
\dspacepl$ and $\ell \leq \ar$. Given any  $\psi \in \dspace$,
a pair of functions $(\phi,\eta) \in \dspace \times \dspace$ is
said to solve the SP on $[\ell(\cdot),\ar(\cdot)]$ for $\psi$
if and only if it satisfies the following properties:
\begin{enumerate}
\item For every $t \in [0,\infty)$, $\phi(t) = \psi(t) +
    \eta(t) \in [\ell(t),r(t)]$;
\item
$\eta = \eta_{\ell} - \eta_{\ar}$, where $\eta_{\ell}$ and $\eta_{\ar}$
are non-decreasing functions such that
\be
\label{comp-sp2}
\int_0^{\infty} \ind_{\{ \phi(s)>\ \ell
(s)\}}d\, \eta_{\ell} (s)=0,\qquad \qquad \int_0^{\infty} \ind_{\{
\phi(s)<\ \ar (s)\}}d\, \eta_{\ar} (s)=0.
\ee
\end{enumerate}
\end{defn}
\noi
If $(\phi,\eta)$ is the unique solution to the SP on $[\ell(\cdot),\ar(\cdot)]$ for
$\psi$ then we will write $\phi = \Gamma_{\ell,\ar} (\psi)$, and
refer to $\Gamma_{\ell,\ar}$ as the associated SM.
Moreover, the pair $(\eta_{\ell}, \eta_{\ar})$ will be referred to as the constraining processes associated
with the SP.

Although Definition \ref{def-smap1} is a natural extension of
the SP to time-dependent domains in $\R$ it is restrictive in
that it only allows ``constraining terms'' $\eta$ that are of
bounded variation. In particular, this implies that any RBM
constructed via the associated SM is automatically a
semimartingale. For fixed domains in $\R^d$, a generalization
of the SP that allows for a pathwise construction of RBMs that
are not necessarily semimartingales was introduced in
\cite{ram1} (see also \cite{burtoby} for a formulation in two
dimensions). The following is the analog of these
generalizations for time-dependent domains in $\R$.

\begin{defn}{\bf (Extended Skorokhod problem on $[\ell(\cdot),r(\cdot)]$)}
\label{def-smap2} Suppose that $\ell  \in \dspacemin$, $\ar \in
\dspacepl$ and $\ell \leq \ar$. Given any  $\psi \in \dspace$,
a pair of functions $(\phi,\eta) \in \dspace \times \dspace$ is
said to solve the ESP on $[\ell(\cdot),\ar(\cdot)]$ for $\psi$
if and only if it satisfies the following properties:
\begin{enumerate}
\item For every $t \in [0,\infty)$, $\phi(t) = \psi(t) +
    \eta(t) \in [\ell(t),r(t)]$;
\item For every $0\leq s\leq t<\infty$, \beas \eta (t) -
    \eta (s) \geq 0 , &&
 \mbox{ if } \phi(u) < r(u) \mbox{ for all } u \in (s,t] \\
 \eta (t) - \eta (s) \leq 0  &&
 \mbox{ if } \phi(u) > \ell(u) \mbox{ for all } u \in (s,t];
 \eeas
\item For every $0\leq t<\infty$, \beas
 \eta (t) - \eta (t-) \geq 0  && \mbox{ if } \phi(t) < r(t) , \\
 \eta (t) - \eta (t-) \leq 0  && \mbox{ if }\phi(t)> \ell(t),
 \eeas
where $\eta(0-)$ is to be interpreted as $0$.
\end{enumerate}
\end{defn}

\noi
If $(\phi,\eta)$ is  the unique solution to the ESP on $[\ell(\cdot),\ar(\cdot)]$ for
$\psi$ then we will write $\phi = \bgamma_{\ell,\ar} (\psi)$, and
refer to $\bgamma_{\ell,\ar}$ as the associated extended Skorokhod map (ESM).  \\

We conclude this section by establishing certain properties of
 SPs and  ESPs (see Theorem 1.3 of \cite{ram1} for analogs for
time-independent multi-dimensional domains).
The first property describes in what sense the ESP is a generalization of the SP.

\begin{prop}
\label{prop:sp-esp}
Suppose we are given $\ell \in \dspacemin$, $\ar \in \dspacepl$ and $\psi \in \dspace$.
If $(\phi, \eta)$ solve the SP on $\gt$ for $\psi$, then $(\phi,\eta)$ solve
the ESP on $\gt$ for $\psi$.   Conversely, if $(\phi, \eta)$ solve the ESP
on $\gt$ for $\psi$ and $\eta$ has finite variation on every bounded interval, then
$(\phi,\eta)$ solve the SP for $\psi$.
\end{prop}
\proof
The first statement follows from the easily verifiable fact that property 2 of
Definition \ref{def-smap1} implies properties 2 and 3 of Definition \ref{def-smap2}.
For the converse, let $(\phi, \eta)$ be a solution to the ESP on $\gt$ for $\psi$ and
suppose $\eta$ has finite variation on every bounded interval.
Then the Lebesgue-Stieltjes measure $d\eta$ is absolutely continuous with respect to the corresponding
total variation measure $d|\eta|$.
Let $\gamma$ be the Radon-Nikod\`{y}m
derivative $d\eta/d|\eta|$ of $d\eta$ with respect to $d|\eta|$. Then
$\gamma$ is $d|\eta|$-measurable, $\gamma(s) \in \{-1, 1\}$ for $d|\eta|$ a.e.
$s\in[0,\infty)$ and
\[
\eta(t) = \int_{[0,t]} \gamma(s)d|\eta|(s). \]
Moreover, it is well-known (see, for example, Section X.4  of \cite{doobbook})
that for $d|\eta|$ a.e.\ $s\in [0,\infty)$,
 \begin{eqnarray} \gamma(s)=\lim_{n\rightarrow \infty}
\frac{\eta(s+\varepsilon_n)-\eta(s-)}{|\eta|(s+\varepsilon_n)-|\eta|(s-)},
\label{gamma}\end{eqnarray} where the sequence $\varepsilon_n$
depends on $s$ and is such that
$|\eta|(s+\varepsilon_n)-|\eta|(s-)>0$ and
$\varepsilon_n\rightarrow 0$ as $n\rightarrow \infty$. Now, for
each $t\geq 0$, define
$\eta_{\ell}(t)=\int_{[0,t]}\ind_{\{\gamma(s)=1\}}d|\eta|(s)$
and
$\eta_{\ar}(t)=\int_{[0,t]}\ind_{\{\gamma(s)=-1\}}d|\eta|(s)$.
Since $\gamma$ only takes the values $1$ and $-1$ ($d|\eta|$
a.e.), it is clear that $\eta=\eta_{\ell}-\eta_{\ar}.$ We shall
now show that $\eta_{\ell}$ satisfies the first complementary
condition in (\ref{comp-sp2}). It follows from the definition
of $\eta_{\ell}$ that
\[ \int_0^{\infty} \ind_{\{ \phi(s)>\ \ell (s)\}}d\, \eta_{\ell} (s)
=\int_0^{\infty} \ind_{\{ \phi(s)>\ \ell
(s)\}}\ind_{\{\gamma(s)=1\}}d|\eta|(s).\] Suppose that there
exists $s\geq 0$ such that $\phi(s)>\ \ell (s)$, $\gamma(s)=1$
and (\ref{gamma}) holds. We will show that this assumption
leads to a contradiction. Since $\phi(s)>\ \ell (s)$, by the
right continuity of $\phi$ and $\ell$, there exists $\delta>0$
such that $\phi(u)>\ \ell (u)$ for all $u\in [s,s+\delta]$. By
properties 2 and 3 of Definition \ref{def-smap2}, we have
$\eta(u)-\eta(s-)\leq 0$ for each $u\in [s,s+\delta]$.  On the
other hand, since $\gamma(s)=1$ and $|\eta|$ is a
non-decreasing function, for all sufficiently large $n$  we
have from (\ref{gamma}) that $\eta(s+\varepsilon_n)-\eta(s-)>
0$. This leads to a contradiction. Hence $\phi(s)>\ \ell (s)$
and $\gamma(s)=1$ cannot hold simultaneously for $d|\eta|$
a.e.\ s, which proves the first complementarity condition in
(\ref{comp-sp2}). The second complementary condition in
(\ref{comp-sp2}) can be established in a similar manner.
\endproof

\begin{cor}
\label{cor-separatelr}
Suppose that $\ell  \in \dspacemin$, $\ar \in
\dspacepl$ and $\inf_{t\geq 0}(\ar(t)-\ell(t))>0$. If $(\phi,\eta) \in
\dspace \times \dspace$ solve the ESP on $[\ell(\cdot),\ar(\cdot)]$ for
some $\psi\in \dspace$, then $(\phi,\eta)$ solve the SP on
$[\ell(\cdot),\ar(\cdot)]$ for $\psi$.
\end{cor}
\proof By Proposition \ref{prop:sp-esp}, it suffices to show
that $\eta$ has bounded variation on every finite time
interval. Let $\tau_0=0$ and for $n\in  \Z_+$, let
$\tau_{2n+1}=\inf\{t\geq \tau_{2n}: \phi(t)=\ar(t)\}$ and
$\tau_{2n+2}=\inf\{t\geq \tau_{2n+1}: \phi(t)=\ell(t)\}$. For
each $n\in \Z_+$, on the interval $[\tau_n,\tau_{n+1})$, $\phi$
will touch exactly one of the boundaries $\ell$ and $\ar$.  By
properties 2 and 3 of the ESP, this implies that $\eta$ will be
either non-decreasing or non-increasing, and hence in
particular of bounded variation, on each interval
$[\tau_n,\tau_{n+1})$.  Moreover, under the assumption
$\inf_{t\geq 0}(\ar(t)-\ell(t))>0$ and the fact that $\phi\in
\dspace$, it is easy to see that  there are finitely many
$\tau_n$'s in each bounded time interval. Thus $\eta$ will have
finite variation on each bounded time interval. \endproof

The following property is a simple, but extremely useful,
closure property of the ESP. Below, the abbreviation u.o.c.\
stands for uniformly on compacts, i.e., we say $f_n \ra f$
u.o.c.\ if for every $T <\infty$, $\sup_{s \in [0,T]} |f_n(s) -
f(s)| \ra 0$ as $n \ra \infty$.

\begin{prop} \label{convergence} {\bf (Closure Property)}
For each $n\in \N$, let $\ell_n  \in \dspacemin$, $\ar_n \in
\dspacepl$ be such that $\ell_n \leq \ar_n$, and let $\psi_n \in
\dspace$.  Suppose there exist $\ell  \in \dspacemin$, $\ar
\in\dspacepl$ and $\psi \in \dspace$ such that $\psi_n \rightarrow
\psi$, $\ell_n\rightarrow \ell$ and $\ar_n\rightarrow \ar$ u.o.c.,
as $n \ra \infty$. Moreover, suppose that for each $n \in \N$,
$(\phi_n,\eta_n)$ solve the ESP on $[\ell_n(\cdot),\ar_n(\cdot)]$
for $\psi_n$.  If $\phi_n \rightarrow \phi$ u.o.c., as
$n\rightarrow \infty$, then $(\phi, \phi - \psi)$ solve the ESP on
$[\ell(\cdot),\ar(\cdot)]$ for $\psi$.
\end{prop}
\proof  Let $\psi_n, \ell_n, \ar_n, \phi_n, \eta_n, n \in \N,$ and $\psi, \ell, \ar, \phi$ be as in
the statement of the proposition and let $\eta \doteq \phi -\psi$.
 By property 1 of Definition \ref{def-smap2},
 $\eta_n = \phi_n - \psi_n$ and $\phi_n(t) \in [\ell_n (t), \ar_n (t)]$
for all $t \in [0,\infty)$.
Together with the assumed u.o.c.\ convergences of $\psi_n, \phi_n, \ell_n$
and $\ar_n$ to $\psi, \phi, \ell$ and $\ar$, respectively, this implies
$\eta_n \ra \eta$ u.o.c., as $n \ra \infty$, and
$\phi(t) \in [\ell (t), \ar(t)]$ for all $t \in [0,\infty)$.
Thus $(\phi,\eta)$ satisfy property 1 of Definition \ref{def-smap2}.

Now, suppose that $\eta(t-) > \eta(t)$ for some $t$. We will show
that then $\phi(t) = r(t)$. Since $\eta_n \to \eta$ u.o.c., we have
$\eta_n(t-) > \eta_n(t)$ for all sufficiently large $n$. By property 3 of Definition
\ref{def-smap2}, this
implies that $\phi_n(t) = r_n(t)$ for all sufficiently large $n$.  The convergences
$\phi_n(t) \to \phi(t)$ and $r_n(t) \to r(t)$ then imply that $\phi(t) =
r(t)$. An analogous argument shows that if $\eta (t-) < \eta (t)$ then
$\phi(t) = \ell(t)$, thus showing that $(\phi,\eta)$ satisfy  property 3 of
Definition \ref{def-smap2}.

In order to show that $(\phi, \eta)$ satisfy the remaining property 2 of Definition \ref{def-smap2},
 fix $0 \leq s < t<\infty$ and suppose that
\begin{equation}
\label{eq-phiare}
\phi(u) < r(u) \qquad \mbox{ for } u\in (s, t].
\end{equation}
 We want to show that then  $\eta(t) \geq \eta(s)$. By the right continuity
of $\eta$, it  suffices to show that $\eta(t) \geq
\eta(\tilde{s})$ for every $\tilde{s} \in (s, t]$. Suppose, to
the contrary, that $\eta (t) < \eta(\tilde{s})$ for some
$\tilde{s} \in (s,t]$. Since $(\phi, \eta)$  satisfy property
3, due to condition (\ref{eq-phiare}) we must have $\eta(u)
\geq \eta(u-)$ for every $u \in [\tilde{s},t]$. In particular,
this implies that the set $B \doteq \{\eta(u): u \in
[\tilde{s}, t]\} \supseteq [\eta(t) , \eta(\tilde{s})]$. Thus
the set  $B$ is uncountable, while  the set $A$ of all
discontinuities of all functions $\ell_n, r_n, \psi_n, \phi_n,
\ell, r, \psi$ and $\phi$ is countable. Hence, there exists
$\alpha \in  [\eta(t),\eta(\tilde{s})) \sm \{\eta(u):u \in
A\}$.  Define $u_1 \doteq \inf \{u \in  (\tilde{s}, t) \sm A:
\eta(u) =\alpha\}$ and note that $u_1 > \tilde{s}$ and
\begin{equation} \label{ineq-closure-eta} \eta(u) > \eta(u_1)
\qquad \mbox{ for all }  u \in [\tilde{s}, u_1).
\end{equation}
In addition, the functions $\phi$ and $\ar$ are continuous at
$u_1$.  Therefore, by (\ref{eq-phiare}) we know that there
exist $u_-\in [\tilde{s}, u_1)$, $u_+ \in (u_1, t]$ such that
$\inf_{u\in [u_-, u_+]} (r(u) - \phi(u)) > 2 \eps$, where $\eps
\doteq (r(u_1) - \phi(u_1))/4 > 0$. Since $\phi_n \to \phi$ and
$r_n \to r$ u.o.c., we know that for all sufficiently large
$n$, $\inf_{u \in [u_-,u_+]} (r_n(u) - \phi_n(u)) > \eps$. Then
property 2 of Definition \ref{def-smap2} implies that
$\eta_n(u_-) \leq \eta_n(u_1)$ for all sufficiently large $n$.
Passing to the limit, we obtain $\eta(u_-) \leq \eta(u_1)$,
which contradicts (\ref{ineq-closure-eta}).  Thus, we must have
$\eta(s) \leq \eta(t)$ when (\ref{eq-phiare}) holds.   An
analogous argument can be used to show that $\eta(s) \geq \eta
(t)$ whenever $\phi(u) > \ell(u)$ for all $u\in [s,t]$. This
completes the proof that $(\phi, \eta)$ solve the ESP on $\gt$
for $\psi$.
\endproof

\subsection{An Explicit Formula for Solutions to the ESP on $\gt$}
\label{subs-explicit}

The following theorem is our main result in this section.

\begin{theorem}
\label{th-sptv} Suppose that $\ell  \in \dspacemin$, $\ar \in
\dspacepl$ and $\ell \leq \ar$. Then for each $\psi\in
\dspace$, there exists a unique pair $(\phi,\eta) \in \dspace
\times \dspace$ that solves the ESP on $\gt$ for $\psi$.
Moreover, the ESM $\bGamma_{\ell, \ar}$ admits the following
explicit representation:
 \be \label{explicit} \bGamma_{\ell,\ar} (\psi)  =
 \psi - \Xi_{\ell,r}(\psi),
  \ee
where the mapping $\Xi_{\ell,r}:\dspace \mapsto \dspace$ is
defined as follows:  for each $t \in [0,\infty)$,
 \be
 \label{def-Xi}
 \ba{l} \Xi_{\ell,r}(\psi)(t) \doteq   \ds \max
 \left( \left[(\psi(0)-r(0))^+\wedge \inf_{u\in
 [0,t]}(\psi(u)-\ell(u))
 \right], \right. \\
 \ds \hspace{1.3in} \left.
 \sup_{s\in [0,t]}\left[(\psi(s)-r(s))\wedge \inf_{u\in
 [s,t]}(\psi(u)-\ell(u))\right] \right).
 \ea
 \ee
Furthermore, the map $(\ell, \ar, \psi) \mapsto \bGamma_{\ell,
\ar}$ is a continuous map on $\dspacemin \times \dspacepl
\times \dspace$ (with respect to the topology of uniform
convergence on compact sets).  Lastly, if $\inf_{t\geq 0}
(\ar(t) - \ell(t)) > 0$ then $\Gamma_{\ell,\ar} =
\bGamma_{\ell, \ar}.$
\end{theorem}

\begin{remark}
\label{rem-explicit}
{\em  When $\ar \equiv \infty$ and $\ell \in \dspace$,
    Definition \ref{def-smap2} reduces to a one-dimensional
    SP with time-varying domain $[\ell(\cdot), \infty)$, and the right-hand side of (\ref{explicit}) reduces to $\Gamma_{\ell} (\psi)$,
where the mapping $\Gamma_{\ell}:\dspace \mapsto \dspace$ is given by
\be
\label{def-gammal}
 \Gamma_\ell (\psi) (t)\doteq \psi (t) + \sup_{s \in [0,t]} \left[\ell(s)
-\psi(s) \right]^+  \mbox{ for } t \in [0,\infty).
\ee
In this situation,  the proof of (\ref{def-g0}) can be
    extended in a straightforward manner (see, for example,
    Lemma 3.1 of \cite{burchesyl04}) to show that $\Gamma_{\ell}$ defines the unique solution
to the associated SP.
}
\end{remark}

The rest of this section is devoted to
the proof of Theorem \ref{th-sptv}.
 For the case of time-independent boundaries
$\ell \equiv 0$ and $\ar \equiv a > 0$, this result was
established in Theorem 2.1  of \cite{krulehramshr2} using a
completely different argument from that used here.
 The proof of Theorem \ref{th-sptv} presented in this paper thus provides, in
particular, an alternative proof of Theorem 2.1 of
\cite{krulehramshr2} (see also \cite{chilaz}, Section 14 of \cite{whibook} and
the discussion in \cite{krulehramshr2} of related formulas in the time-independent case).

For the rest of this section, we fix $\ell  \in \dspacemin$ and
$\ar \in \dspacepl$ such that $\ell \leq \ar$.
We first establish uniqueness of solutions to the ESP on
$\gt$ in Proposition \ref{prop-spuniq} ---the proof is a relatively
straightforward modification of the standard proof for the SP
on $[0,\infty)$ (see, for example, Lemma 3.6.14 in
\cite{karshrbook} and also Lemma 3.1  of \cite{burchesyl04}).

\begin{prop}
\label{prop-spuniq} Given any $\psi \in \dspace$, there exists
at most one $\phi \in \dspace$ that satisfies the ESP on $\gt$
for $\psi$.
\end{prop}
\begin{proof}
Let $(\phi,\eta)$ and $(\phi',\eta')$ be two pairs of functions in
$\dspace \times \dspace$ that solve the ESP on
$[\ell(\cdot),r(\cdot)]$ for $\psi\in \dspace$. Suppose that there
exists $T\geq 0$ such that $\phi(T)>\phi'(T)$. Let \be
\label{def1-tau} \tau=\sup\{ t\in[0, T]:\ \phi(t)\leq \phi'(t)\}.
\ee Then it follows that $\phi(\tau-)\leq \phi'(\tau-)$.  We now
consider two cases.

\noindent {\bf Case 1.} $\phi(\tau)\leq \phi'(\tau)$. In this
case, for $t\in (\tau,T]$, by the definition of $\tau$ and
property 1 of Definition \ref{def-smap2}, we have $\ell(t)\leq
\phi'(t)<\phi(t)\leq r(t)$. Since on $(\tau,T]$, $\phi$ will not
hit $\ell$ and $\bphi'$ will not hit $\ar$, by property 2 of
Definition \ref{def-smap2}, we see that $\eta(T)-\eta(\tau)\leq 0$
and $\eta'(T)-\eta'(\tau)\geq 0$. Consequently,
 \[
0<\phi(T)-\phi'(T)=
\eta(T)-\eta'(T) \leq  \eta(\tau)-\eta'(\tau) = \phi(\tau)-\phi'(\tau),
\]
which contradicts the case assumption.

\noindent {\bf Case 2.} $\phi(\tau)> \phi'(\tau)$. In this
case,  we have $\phi(\tau) > \ell (\tau)$ and $\phi'(\tau) < \ar (\tau)$.
By property 3 of Definition \ref{def-smap2}, this implies that
$\eta(\tau) - \eta(\tau-) \leq 0$ and $\eta'(\tau) - \eta'(\tau-) \geq 0$.
When combined with property 1 of Definition \ref{def-smap2}, this shows that
\[ 0 < \phi (\tau) - \phi'(\tau) = \eta(\tau) - \eta' (\tau) \leq
\eta (\tau-) - \eta' (\tau-) = \phi(\tau-) - \phi' (\tau-),
\]
which contradicts the definition (\ref{def1-tau}) of $\tau$.

We thus conclude that $\phi(T)\leq \phi'(T)$ for all $T \geq 0$.
Using an exactly analogous argument we can show that $\phi'(T)\leq \phi(T)$ for
all $T\geq 0$. Hence $\phi(T)= \phi'(T)$ and, therefore,
$\eta(T)= \eta'(T)$ for all $T \geq 0$.
\end{proof}

Next, in Proposition \ref{piececonstant}, we show that the ESM  is given by the formula (\ref{explicit})
when $\ell, \ar$ and $\psi$
are piecewise constant. The proof will make use of the following family of mappings:
given $\ell \in \dspacemin, \ar \in \dspacepl$ with $\ell \leq \ar$,
for $t \in [0,\infty)$, consider the mapping $\pi_t:\R \ra \R$ with the property that
$\pi_t(x) = x$ if $x \in [\ell(t),\ar(t)]$,  $\pi_t(x) \in \{\ell (t),
  \ar(t)\}$ if $x \not \in [\ell(t), \ar(t)]$ and
\be
\label{def-proj}
\begin{array}{rl}
\pi_t(x) - x \geq 0 & \mbox{ if }  \pi_t (x) = \ell (t) , \\
\pi_t(x) - x \leq 0 & \mbox{ if } \pi_t(x) = \ar(t).
\end{array}
\ee It is straightforward to deduce that, for every $t \geq 0$,
there exists a unique mapping with these properties that is given
explicitly by
 \be \label{implication-baretazero}
 \pi_t(x) = x + [\ell(t)-x]^+-[x-r(t)]^+ = (x \wedge \ar(t) ) \vee \ell(t).
 \ee
Using property 3 of the ESP, it is easy to verify that the ESM
$\bGamma_{\ell,\ar}$ must satisfy \be \label{jump}
\bGamma_{\ell,\ar}(\psi)(0) = \pi_0 (\psi(0)), \qquad
\bGamma_{\ell, \ar} (\phi)(t) = \pi_t \left( \bGamma_{\ell,
\ar}(\psi)(t-) + \psi(t) - \psi(t-) \right)  \quad \forall t >
0. \ee

\begin{prop} \label{piececonstant}
Suppose that $\ell$, $\ar$ and $\psi$ are three piecewise constant
functions in $\dspacemin$, $\dspacepl$ and $\dspace$,
respectively, each with a finite number of jumps and such that $\ell \leq \ar$.
Then for each $\psi\in \dspace$, the pair
$(\psi-\Xi_{\ell,\ar}(\psi),-\Xi_{\ell,\ar}(\psi))$ is the unique
solution to the ESP on $[\ell(\cdot),\ar(\cdot)]$, i.e.,
$\bGamma_{\ell,\ar}(\psi)=\psi-\Xi_{\ell,\ar}(\psi)$.
\end{prop}
\proof Fix $\psi, \ell, \ar$ as in the statement of the
proposition, and let $\{\pi_t, t \in [0,\infty)\}$ be the
associated family of mappings as defined in
(\ref{implication-baretazero}). Now, let $J= \{t_1, t_2,
\ldots, t_n \}$ be the union of the times of jumps of $\ell$,
$\ar$ and $\psi$, suppose $0<t_1<t_2<\cdots <t_n<\infty$, and
set $t_{n+1} \doteq \infty$. Define $\phi \doteq \psi -
\Xi_{\ell, \ar} (\psi)$ and $\eta \doteq \phi -\psi$. We will
use induction to show that $(\phi, \eta)$ solve the ESP on
$\gt$ for $\psi$. When $t =0$, it is straightforward to verify
from (\ref{implication-baretazero}) and the definition of
$\Xi_{\ell, \ar}$ that $\phi (0) = (\psi(0) \wedge \ar(0)) \vee
\ell (0) = \pi_0 (\psi(0))$. When combined with (\ref{jump}),
this shows that $(\phi, \eta)$ solve the ESP (on $[\ell
(\cdot), \ar (\cdot)]$) for $\psi$  when $t = 0$. Since $\ell,
\ar, \psi$ are constant on $[0,t_{1})$, it immediately follows
from the definition (\ref{def-Xi}) of $\Xi_{\ell,\ar}$ that
$\phi$ is also constant on $[0,t_{1})$, and so it follows that
$(\phi,\eta)$ solve the ESP for $\psi$ on $[0,t_{1})$.

Now, suppose $(\phi, \eta)$ solve the ESP  on
$[\ell(\cdot),\ar(\cdot)]$ for $\psi$ over the time interval $[0,
t_m)$ for some $m \in \{1, \ldots, n\}$. We first observe that, for any $t \in
[0,\infty)$, $\Xi_{\ell,\ar}(\psi)(t)$ is the maximum of the following three terms:
\begin{enumerate}
\item
$(\psi(0)-\ar(0))^+\wedge \inf_{u\in [0,t)}(\psi(u)-\ell(u))\wedge (\psi(t)-\ell(t))$
\item
$\sup_{s\in [0,t)}\left[(\psi(s)-\ar(s))\wedge \inf_{u\in [s,t)}(\psi(u)-\ell(u))\wedge (\psi(t)-\ell(t))\right]$,
\item
$(\psi(t)-\ar(t))\wedge (\psi(t)-\ell(t))$
\end{enumerate}
and therefore admits the representation
\[\Xi_{\ell,\ar}(\psi)(t) =  \max\left[ \Xi_{\ell, \ar}(\psi) (t-), (\psi
(t) - \ar(t))\right] \wedge (\psi(t)-\ell(t)).
\]
Recalling the description of the map $\pi_t$ given in  (\ref{implication-baretazero}),
we see that
\begin{eqnarray*}
 \phi(t) & = &  \psi(t) - \max\left( \Xi_{\ell, \ar}(\psi) (t-), (\psi (t) - \ar(t))\right)\wedge (\psi(t)-\ell(t)) \\
& = & \min\left(\psi(t) - \Xi_{\ell, \ar} (\psi) (t-), \ar(t) \right) \vee \ell (t) \\
& = & \pi_t (\psi(t) - \Xi_{\ell, \ar} (\psi) (t-)) \\
& = & \pi_t (\phi(t-) + \psi(t) -\psi(t-)).
\end{eqnarray*}
Substituting  $t = t_m$, this yields the relation $\phi (t_m) = \pi_{t_m}
(\phi(t_m-) + \psi(t_m) - \psi(t_m-))$.  By (\ref{jump}), this implies
$(\phi, \eta)$ solve the ESP on $[\ell(\cdot), \ar(\cdot)]$ for $\psi$ during the
interval $[0,t_m]$.
Once again, since $\psi, \ell, \ar$, and therefore $\phi$, are
constant on $[t_m, t_{m+1})$ this implies that $(\phi, \eta)$
solve the ESP on $\gt$ for $\psi$ on $[0,t_{m+1})$. By the
induction argument and the uniqueness result established in
Proposition \ref{prop-spuniq}, we have the desired
result.\endproof

 A simple approximation argument can now be used
 to complete the proof of Theorem \ref{th-sptv}.

\begin{proof}[Proof of Theorem \ref{th-sptv}]
Given $\ell \in \dspacemin, \ar\in \dspacepl$ such that $\ell
\leq \ar$, it is easy to see that there exist sequences of
functions  $\ell_n \in \dspacemin$, $n \in \N$, $\ar_n \in
\dspacepl$, $n \in \N$, with $\ell_n \leq \ar_n$, that are
piecewise constant with a finite number of jumps and such that
$\ell_n\rightarrow \ell$, $\ar_n\rightarrow \ar$ u.o.c.\  as
$n\rightarrow \infty$. Likewise, given $\psi_n \in \dspace$,
there exists a sequence of piecewise constant functions
$\psi_n$ with a finite number of jumps such that $\psi_n \ra
\psi$ u.o.c., as $n \ra \infty$. For each $n\in  \N$, by
Proposition \ref{piececonstant}, we know that
$\Gamma_{\ell_n,\ar_n}(\psi_n)= \phi_n \doteq
\psi_n-\Xi_{\ell_n,\ar_n}(\psi_n)$. Since $\psi_n-\ell_n$ and
$\psi_n-\ar_n$ converge u.o.c., as $n\rightarrow \infty$, to
$\psi-\ell$ and $\psi-\ar$, respectively, and u.o.c.
convergence is preserved under the operations $\inf,\ \sup,\
\wedge,\ \max$, we then conclude, from (\ref{def-Xi}), that
$\phi_n = \psi_n - \Xi_{\ell_n, \ar_n}(\psi_n)\rightarrow \psi
- \Xi_{\ell, \ar}(\psi)$ u.o.c, as $n \ra \infty$. In
particular, it is clear that $\Xi_{\ell, \ar}$ is a continuous
map on $\dspace$ (with respect to the topology of u.o.c.
convergence). By the closure property (Proposition
\ref{convergence}),
$(\psi-\Xi_{\ell,\ar}(\psi),-\Xi_{\ell,\ar}(\psi))$ is a
solution to the ESP on $[\ell(\cdot),\ar(\cdot)]$ for $\psi$.
Uniqueness follows from  Proposition \ref{prop-spuniq}. In
particular, this shows that the map $(\ell, \ar, \psi)\mapsto
\bGamma_{\ell, \ar}(\psi)$  is continuous with respect to the
topology of u.o.c. convergence. The last assertion of the
theorem is a direct consequence of Corollary
\ref{cor-separatelr}.
\end{proof}

\section{Comparison results}
\label{sec-comp}

This section presents some ``comparison'' or ``monotonicity''
results. They are quite intuitive but their proofs require some
technical arguments. Recall the definition of the pair of constraining
processes $(\eta_\ell, \eta_\ar)$ associated with an SP given in Definition \ref{def-smap1}.
Section \ref{subs-mondom} establishes
monotonicity of the individual constraining processes with
respect to the domain $\gt$ for a fixed $\psi$, while in Section
\ref{subs-monpsi}, monotonicity of the constraining
processes with respect to the input $\psi$ is established for a given
time-varying domain $\gt$.

\subsection{Monotonicity with respect to the domain}
\label{subs-mondom}

The main result, Proposition \ref{prop-monot}, will be preceded by  a few lemmas.

\begin{lemma}
 \label{comparison1}
Assume that $\ell,\tilde{\ell}  \in \dspacemin$,
$\ar,\tilde{\ar} \in \dspacepl$, and $\tilde{\ell} = \ell$, $\ar
\leq \tilde{\ar}$ and $\inf_{t\geq 0} (\ar(t) - \ell(t)) > 0$.
Let $\Gamma_{\ell,\ar}$ and $\Gamma_{\tilde
\ell,\tilde \ar}$ be the associated SMs on
$[\ell(\cdot),\ar(\cdot)]$ and $[\tilde \ell(\cdot),\tilde
\ar(\cdot)]$, respectively, and given $\psi \in \dspace$, let
$(\eta_{\ell}, \eta_{\ar})$ and $({\eta}_{\tilde{\ell}},
{\eta}_{\tilde{\ar}})$ be the corresponding pairs of
constraining processes.  Then, for every $t \in [0,\infty)$,
\be
\label{comp-decomp}
  \eta_{\ar} (t)\geq  \eta_{\tilde \ar}(t) \quad \quad \mbox{  and  }
 \quad \quad  \eta_{\ell} (t)\geq \eta_{\tilde \ell}(t).
\ee
\end{lemma}

\begin{proof} Let $ \tilde \ell , \ell , \ar , \tilde{\ar} $ and $\psi$
be as in the statement of the lemma. First note that by Theorem
\ref{th-sptv}, the conditions on $\tilde{\ell}, \tilde{\ar},
\ell, \ar$ guarantee that solutions to the SP  on both $[\ell,
\ar]$ and $[\tilde{\ell}, \tilde{\ar}]$ exist for all $\psi \in
\dspace$  and so the pairs of constraining processes
$(\eta_{\ell}, \eta_{\ar})$ and $(\eta_{\tilde{\ell}},
\eta_{\tilde{\ar}})$ are well-defined. Moreover,  the explicit
formula (\ref{explicit})  of Theorem \ref{th-sptv}, when
combined with the decomposition $\eta= \Gamma_{\ell, \ar}
(\psi) -\psi = \eta_{\ell} -\eta_{\ar}$ (see Definition
\ref{def-smap1}), shows that
\begin{equation}
\label{eq-intermed}
 \eta_\ar = \eta_\ell + \Xi_{\ell,\ar}(\psi) \quad \quad \mbox{ and } \quad \quad
 \eta_{\tilde{r}} = \eta_{\tilde{\ell}} + \Xi_{\tilde{\ell},\tilde{\ar}} (\psi),
\end{equation}
where $\Xi$ is as defined in (\ref{def-Xi}).
For a fixed $\ell$, it is easily verified
from the explicit formula (\ref{def-Xi}) that the map
$\ar \mapsto \Xi_{\ell,\ar}(\psi)$ is monotone non-increasing (with respect to the obvious ordering).
 Since $\ell = \tilde{\ell}$ and $\ar \leq \tilde{\ar}$, this implies
\be
\label{comp-xi}
\Xi_{\tilde \ell,\tilde \ar}(\psi)\leq
\Xi_{\ell,\ar}(\psi) \qquad \forall  \psi \in \dspace.
\ee
By (\ref{explicit}), this is equivalent to the relation
 \be \label{ineq-xi} \Gamma_{\tilde \ell,\tilde \ar}(\psi) \geq
\Gamma_{\ell,\ar}(\psi) \qquad \forall  \psi \in \dspace.
 \ee
Combining (\ref{eq-intermed}) and (\ref{comp-xi}), it follows
that in order to show (\ref{comp-decomp}), it suffices to show
that
\begin{equation}
\label{comp-toshow}
 \eta_\ell(t)\geq \eta_{\tilde \ell}(t)  \quad \forall \, t\geq 0.
\end{equation}

Since $\Gamma_{\ell, \ar}(\psi) (0) = \pi_0 (0)$ by
(\ref{jump}) and $\ell(0) < \ar(0)$, from the complementarity
conditions (\ref{comp-sp2}) it is clear that $\eta_{\ell} (0) =
[ \ell(0) - \psi(0)]^+$ with the analogous expressions for
$\eta_{\tilde \ell}$. Since  $\ell = \tilde{\ell}$, this
immediately implies (\ref{comp-toshow}), in fact with equality,
for $t = 0$.

Now, let
\[t^* \doteq \inf\{s\geq 0:  \eta_\ell(s) <  \eta_{\tilde \ell}(s)\}. \]
We will argue by contradiction to show that $t^* = \infty$.
Indeed, suppose that $t^*<\infty$. Then
 \be \label{contra1a}
 \eta_\ell(t^*-)\geq  \eta_{\tilde \ell}(t^*-)
\ee
and for all $\ve_0 >0$ there exists  $\eps\in(0,\ve_0)$ such that
\be
\label{contra1b}
\eta_\ell(t^*+\epsilon)< \eta_{\tilde \ell}(t^*+\epsilon).
 \ee
Invoking (\ref{jump}) and the inequality
  $\Gamma_{\tilde \ell,\tilde \ar}(\psi)(t^*-)\geq
\Gamma_{\ell,\ar}(\psi)(t^*-)$ from (\ref{ineq-xi}),
we obtain
\[
\ba{rcl}  \eta_\ell(t^*)- \eta_\ell(t^*-) & = &
[\ell(t^*)-\Gamma_{\ell,\ar}
(\psi) (t^*-)-(\psi(t^*)-\psi(t^*-))]^+ \\
& \geq  & [\tilde{\ell}(t^*)-\Gamma_{\tilde \ell,\tilde \ar}
(\psi) (t^*-)-(\psi(t^*)-\psi(t^*-))]^+ \\
& = & \eta_{\tilde \ell}(t^*)- \eta_{\tilde
\ell}(t^*-).
\ea
\]
When combined with (\ref{contra1a}), this
implies that $ \eta_\ell(t^*)\geq  \eta_{\tilde \ell}(t^*)$. We
now consider two cases.  If $ \eta_\ell(t^*)>  \eta_{\tilde
\ell}(t^*)$, then  the right-continuity of $ \eta_\ell$ and $
\eta_{\tilde \ell}$ dictates that $ \eta_\ell(t^*+\epsilon)
>  \eta_{\tilde \ell}(t^*+\varepsilon)$ for every positive
$\varepsilon$ small enough, which contradicts (\ref{contra1b}).
On the other hand,  suppose $ \eta_\ell(t^*)= \eta_{\tilde
\ell}(t^*)$.  When combined with (\ref{contra1b}) and the fact
that $\eta_{\ell}$ is non-decreasing, this implies that for
every $\ve_0 > 0$, there exists $\ve \in (0,\ve_0)$ such that
$\eta_{\tilde \ell}(t^*) < \eta_{\tilde \ell}(t^* +\ve)$. Due
to the complementarity condition (\ref{comp-sp2}) and the
right-continuity of $\Gamma_{\tilde \ell, \tilde \ar} (\psi)$,
this, in turn,  implies that $\Gamma_{\tilde \ell,\tilde \ar}
(\psi) (t^*)= \tilde{\ell} (t^*)= \ell (t^*)$. Since
$\Gamma_{\tilde \ell,\tilde \ar}(\psi) \geq
\Gamma_{\ell,\ar}(\psi) \geq \ell$, this means that
$\Gamma_{\ell,\ar}(\psi)(t^*) =\Gamma_{\tilde \ell,\tilde
\ar}(\psi)(t^*)=\ell(t^*)$. Along with  the relation $\tilde
\ell(t^*) \leq \ell(t^*) < \ar (t^*) \leq \tilde{\ar} (t^*)$,
the right-continuity of $\ell$ and $\ar$ and the definition of
the SP, it is easy to see that this implies that for all
sufficiently small $\ve$, $\Gamma_{\tilde \ell,\tilde
\ar}(\psi)(t^*+\ve)$ (respectively, $\Gamma_{ \ell, \ar}(\psi)
(t^* + \ve)$) is equal to
 $\Gamma_{\tilde \ell} (\psi^{*})(\ve)$ (respectively,
$\Gamma_\ell (\psi^*)(\ve)$),
where $\Gamma_{\ell}$ is as defined in (\ref{def-gammal}) and
\[\psi^{*} (t) = \ell (t^*) + \psi(t^*+ t) - \psi(t^*) \quad \quad \mbox{ for } t \geq 0. \]
In particular, using (\ref{def-gammal}), this shows that for all  $\ve$ sufficiently
small,
\[
\ba{rcl}
 \eta_{\tilde \ell} (t^* + \ve) -  \eta_{\tilde \ell} (t^*) & = &
 \sup_{s \in [0,\ve]} \left[\tilde \ell(s) -  \ell (t^*) -
\psi(t^*+s) + \psi(t^*) \right] \\
& \leq &
 \sup_{s \in [0,\ve]} \left[\ell(s) -
\ell(t^*) -
\psi(t^*+s) + \psi(t^*) \right] \\
& = & \eta_{\ell} (t^* + \ve) -  \eta_{\ell} (t^*).
\ea
\]
Since we are considering the case $\eta_{\ell} (t^*) =\eta_{\tilde \ell} (t^*)$,
this once again contradicts
 (\ref{contra1b}). Thus we have shown that
$t^*=\infty$, and hence that (\ref{comp-toshow}) holds.
\end{proof}

\begin{cor}
 \label{comparison2}
Assume that $\ell,\tilde{\ell}  \in \dspacemin$,
$\ar,\tilde{\ar} \in \dspacepl$,  $\tilde{\ell} \leq \ell$,
$\ar = \tilde{\ar}$ and $\inf_{t \geq 0} (\ar(t) - \ell(t)) > 0$.
Let $\Gamma_{\ell,\ar}$ and $\Gamma_{\tilde
\ell,\tilde \ar}$ be the associated SMs on
$[\ell(\cdot),\ar(\cdot)]$ and $[\tilde \ell(\cdot),\tilde
\ar(\cdot)]$, respectively, and given $\psi \in \dspace$, let
$(\eta_{\ell}, \eta_{\ar})$ and $({\eta}_{\tilde{\ell}},
{\eta}_{\tilde{\ar}})$ be the corresponding pairs of
constraining processes.  Then, for every $t \in [0,\infty)$,
\[  \eta_{\ar} (t)\geq  \eta_{\tilde \ar}(t) \quad \quad \mbox{  and  }
 \quad \quad  \eta_{\ell} (t)\geq \eta_{\tilde \ell}(t).
\]
\end{cor}

\begin{proof}
The corollary follows from Lemma \ref{comparison1} by
multiplying all functions by $-1$.
\end{proof}

\begin{prop}
\label{prop-monot} Assume that $\ell,\tilde{\ell}  \in
\dspacemin$, $\ar,\tilde{\ar} \in \dspacepl$, $\tilde{\ell}
\leq \ell$, $\ar \leq \tilde{\ar}$ and $\inf_{t \geq 0} (\ar(t)
- \ell(t)) > 0$. Let $\Gamma_{\ell,\ar}$ and $\Gamma_{\tilde
\ell,\tilde \ar}$ be the associated SMs on
$[\ell(\cdot),\ar(\cdot)]$ and $[\tilde \ell(\cdot),\tilde
\ar(\cdot)]$, respectively, and given $\psi \in \dspace$, let
$(\eta_{\ell}, \eta_{\ar})$ and $({\eta}_{\tilde{\ell}},
{\eta}_{\tilde{\ar}})$ be the corresponding pairs of
constraining processes.  Then, for every $t \in [0,\infty)$,
\[  \eta_{\ar} (t)\geq  \eta_{\tilde \ar}(t) \quad \quad \mbox{  and  }
 \quad \quad  \eta_{\ell} (t)\geq \eta_{\tilde \ell}(t).
\]
\end{prop}
\begin{proof}
Let $\Gamma_{\ell, \tilde \ar}$ be the SM on
$[\ell(\cdot),\tilde \ar(\cdot)]$, and given $\psi \in
\dspace$, let $(\eta_{\ell}^*, \eta_{\tilde \ar}^*)$ be the
corresponding vector of constraining processes. Since $r \leq \tilde{r}$ and $\inf_{t\geq 0} (\tilde{r}(t) - \ell(t)) >  0$
by Lemma \ref{comparison1}, we know that
 \[  \eta_{\ar} (t)\geq
\eta_{\tilde \ar}^*(t) \quad \quad \mbox{  and  }
 \quad \quad  \eta_{\ell} (t)\geq
\eta_{\ell}^*(t).
 \]
Similarly, when $\Gamma_{\ell, \tilde \ar}$ and $\Gamma_{\tilde \ell,
\tilde \ar}$ are considered, by Corollary \ref{comparison2} we obtain
 \[
 \eta_{\tilde \ar}^* (t)\geq  \eta_{\tilde \ar}(t) \quad
\quad \mbox{  and  }
 \quad \quad  \eta_{\ell}^* (t)\geq
\eta_{\tilde \ell}(t).
\]
When combined, these inequalities yield the desired result.
\end{proof}

\subsection{Monotonicity with respect to input trajectories}
\label{subs-monpsi}

Given a fixed time-dependent domain $\gt$, in Proposition \ref{lem:Mono} we
first establish the monotonicity of  $\bGamma_{\ell, \ar} (\psi)$  and the net
constraining term $\bGamma_{\ell, \ar} (\psi ) - \psi$
with respect to input trajectories $\psi$.
For the case when  $\gt = [0,a]$ for some $a > 0$,
this result was established as Theorem 1.7 of
\cite{krulehramshr1}.   Here, we use a simpler argument involving
approximations  to prove the more general result.

\begin{prop} \label{lem:Mono}
Given $\ell,r\in \dspace$ with $\ell \leq \ar$, $c_0, c_0' \in \R$ and
$\psi,\psi'\in \dspace$, suppose
$(\phi,\eta)$ and $(\phi',\eta')$ solve the ESP
 on $[\ell(\cdot),r(\cdot)]$ for $c_0 + \psi$ and $c_0'+ \psi'$,
respectively. If $\psi=\psi'+\nu$ for some non-decreasing function
$\nu\in \dspace$ with $\nu(0)=0$,
 then for each $t\geq 0$, the following two relations hold:
\begin{enumerate}
\item
$[-[c_0-c_0']^+-\nu(t)]\vee [-(\ar(t)-\ell(t))] \leq
    \phi'(t)-\phi(t)\leq [c_0'-c_0]^+ \wedge
    [\ar(t)-\ell(t)];$
\item $\eta(t) -[c_0'-c_0]^+\leq \eta'(t)\leq \eta(t) +
    \nu+[c_0-c_0']^+. $
\end{enumerate}
\end{prop}
\begin{proof}
We first establish property 1 under the additional assumption
that the functions $\ell, \ar, \psi, \psi^\prime$ and $\nu$ stated in the lemma
are piecewise constant with a finite number of jumps.
Since $\phi(t), \phi^\prime(t)$ lie in $[\ell(t), \ar(t)]$ for every $t\in [0,\infty)$,
in order to show the first property it suffices to show that
\be
 -[c_0-c_0']^+-\nu(t)\leq
\phi'(t)-\phi(t)\leq [c_0'-c_0]^+.
\label{ineq:1}
\ee
Let $t_0 =0$ and let $t_1 <  t_2 <  \ldots <
t_{m}$ be the ordered jump times of all the functions
$\ell, \ar, \psi, \psi^\prime$ and $\nu$.
Recall the family of  (time-dependent) projection operators
 $\pi_t$, $t \geq 0$, defined in (\ref{implication-baretazero}).
Using the explicit expression for $\pi_t$, a simple case-by-case verification
shows that for every $t \geq 0$ and
$x, y \in \R$,
\begin{equation}
\label{ineq-key}
 -[y-x]^+\leq \pi_t(x) - \pi_t(y)  \leq [x-y]^+.
\end{equation}
By (\ref{jump}) and the piecewise constant nature of the functions, it follows that
$\phi(t) = \phi(0) =\pi_0(c_0)$ and  $\phi'(t) = \phi'(0) = \pi_0 (c_0')$ for $t \in [0,t_1)$.
When combined with (\ref{ineq-key}), this shows that (\ref{ineq:1}) holds for $t \in [0,t_1)$.
Now suppose that (\ref{ineq:1}) holds for $t\in [0,t_{k-1})$ for some $k \in \{2, \ldots, m\}$.
Then, by  (\ref{jump}) we know that for $t \in [t_{k-1}, t_{k})$,
\begin{equation}
\label{comp-pi}
  \phi (t) = \pi_{t_k} \left( \phi\left(t_{k-1}\right)
+ \psi\left(t_k\right)
- \psi\left(t_{k-1}\right)  \right), \ \
  \phi' (t) = \pi_{t_k} \left( \phi'\left(t_{k-1}\right)
+ \psi'\left(t_k\right)
- \psi'\left(t_{k-1}\right)  \right).
\end{equation}
Another application of (\ref{ineq-key}), along with the relation
$- [z]^+ = [-z] \wedge 0$ and the
fact that $\psi = \psi' + \nu$, implies that
\[ 0 \wedge [\phi'(t_{k-1}) - \phi(t_{k-1}) + \nu(t_{k-1}) - \nu(t_k)]  \leq \phi' (t_k) - \phi(t_k) \leq
 [\phi' (t_{k-1}) - \phi (t_{k-1})  + \nu(t_{k-1}) - \nu(t_{k}) ]^+.
\]
The function $\nu$ is non-decreasing and non-negative and
by the induction assumption, the first inequality in (\ref{ineq:1}) holds for $t = t_{k-1}$.
Therefore
\[ -[c_0 - c_0']^+ - \nu (t_k) \leq \phi'(t_k) - \phi (t_k) \leq [c_0' - c_0]^+. \]
Since $\phi$ and $\eta$ are constant on $[t_k, t_{k+1})$, we have shown
that (\ref{ineq:1}) holds for $t \in [0,t_{k+1})$ and, by induction, for $t \in [0,\infty)$
when all the relevant functions are piecewise constant.

  For the general case, let $\ell_n \in \dspacemin, \ar_n \in \dspacepl$, $n \in \N$, be
sequences of piecewise constant functions with a finite number
of jumps such that  $\ell_n \leq \ar_n$ for every $n \in \N$
and $\ell_n \ra \ell$ and $\ar_n \ra \ar$ u.o.c., as $n \ra
\infty$. Moreover,  let $\psi_n, \psi^\prime_n, \nu_n \in
\dspace$, $n \in \N$, be sequences of piecewise constant
functions with a finite number of jumps such that $\nu_n$ is
non-decreasing and $\psi_n =\psi^\prime_n +\nu_n$  and $\psi_n
\ra \psi$, $\psi_n^\prime \ra \psi^\prime$ u.o.c., as $n \ra
\infty$ (see the proof of Lemma 3.3 in \cite{ramrei2} for an
explicit construction that shows such sequences exist).
Moreover, let $\phi_n = \bGamma_{\ell, \ar} (c_0 + \psi_n)$ and
$\phi'_n = \bGamma_{\ell, \ar}(c_0' + \psi'_n)$. Then the
continuity of the map $\psi \mapsto \bGamma_{\ell, \ar} (\psi)$
established in Theorem \ref{th-sptv} shows that $\phi_n \ra
\phi$ and $\phi'_n \ra \phi'$ u.o.c., as $n \ra \infty$.
Furthermore,  the arguments in the previous paragraph show that
for every $n \in \N$ and $t \in [0,\infty)$, (\ref{ineq:1})
holds with $\phi, \eta$ replaced by $\phi_n$ and $\nu_n$,
respectively. Taking limits as $n \ra \infty$, we obtain
property 1.

The second property can be deduced from the first  using the basic relation
\[\eta'-\eta=\phi'-\phi-(c_0'-c_0)-(\psi'-\psi)=\phi'-\phi-(c_0'-c_0)+\nu.
\]
\end{proof}

Next, we establish monotonicity of the individual constraining
processes $\eta_{\ell}$ and $\eta_{\ar}$ with respect to input
trajectories $\psi$.

\begin{prop}
\label{comp} Given $\ell \in \dspacemin$, $\ar \in \dspacepl$
satisfying $\inf_{t\geq 0} (\ar(t) - \ell(t)) > 0$,
$c_0,c_0'\in \R$ and $\psi,\psi'\in \dspace$ with
$\psi(0)=\psi'(0)$, suppose $(\phi,\eta)$ and $(\phi',\eta')$
solve the SP on $[\ell(\cdot),\ar(\cdot)]$ for $c_0+\psi$ and
$c_0'+\psi'$, respectively. Moreover, suppose $(\eta_\ell,
\eta_\ar)$ and $(\eta'_\ell, \eta'_\ar)$ are the corresponding
constraining processes. If there exists a non-decreasing
function $\nu$ with $\nu(0)=0$ such that $\psi=\psi'+\nu$, then
for each $t\geq 0$, the following two relations hold:
\begin{enumerate}
\item $\eta_\ell(t) -[c_0'-c_0]^+\leq \eta_\ell'(t)\leq
    \eta_\ell(t) + \nu(t) +[c_0-c_0']^+; $
\item $\eta_r'(t) -[c_0'-c_0]^+\leq \eta_r(t)\leq
    \eta_r'(t) + \nu(t)+[c_0-c_0']^+. $
\end{enumerate}
\end{prop}

\begin{proof} Fix $t\in [0,\infty)$.
 Define
\[ \alpha\doteq\inf\{t> 0: \eta_\ell(t)+\nu(t)+[c_0-c_0']^+< \eta_\ell'(t)
\mbox{ or } \eta_\ar(t)+[c_0'-c_0]^+<\eta_\ar'(t) \}, \] where $\alpha=\infty$ if the infimum is over the empty set. Then it follows that for each $s\in [0,\alpha)$, the following two inequalities hold:
\begin{eqnarray}
\eta_\ell'(s)&\leq &\eta_\ell(s)+\nu(s)+[c_0-c_0']^+, \label{est:20}\\
\eta'_\ar(s) & \leq & \eta_\ar(s)+[c_0'-c_0]^+. \label{est:1}
\end{eqnarray}
Suppose $\alpha<\infty$. Then we claim and prove below that the following relations are satisfied:
\begin{eqnarray}
\eta_\ell'(\alpha)&\leq &\eta_\ell(\alpha)+\nu(\alpha)+[c_0-c_0']^+, \label{est:4}\\
\eta'_\ar(\alpha) & \leq & \eta_\ar(\alpha)+[c_0'-c_0]^+. \label{est:5}
\end{eqnarray}
It is easy to see from (\ref{est:20}), (\ref{est:1}) and the
non-decreasing property of $\eta_\ell$ and $\eta_\ar$ that if
$\eta_{\ell}'$ (respectively, $\eta_\ar'$) is continuous, then
(\ref{est:4}) (respectively, (\ref{est:5})) holds. Thus the
claim holds if both $\eta_{\ell}'$ and $\eta_{\ar}'$ are
continuous. We now prove the claim under the assumption that
$\eta_\ell'(\alpha)-\eta_\ell'(\alpha-)>0$. {F}irst note that
by the complementarity condition in (\ref{comp-sp2}) we have
$\phi'(\alpha)=\ell(\alpha)$ and $\eta_\ar'$ is continuous at
$\alpha$, and hence (\ref{est:5}) holds. It follows that
\begin{eqnarray}\eta'_\ell(\alpha)&=&
\eta_\ell'(\alpha-)+\psi'(\alpha-)-\phi'(\alpha-)-\psi'(\alpha)+\ell(\alpha)
\nonumber\\ & =&
-c_0'+\eta_\ar'(\alpha-)-\psi'(\alpha)+\ell(\alpha) \nonumber\\
&=&-c_0'+\eta_\ar'(\alpha-)-\psi(\alpha)+\nu(\alpha)+\ell(\alpha).
\label{est:21}\end{eqnarray} Since
$\eta_\ar(\alpha-)=c_0+\psi(\alpha-)+\eta_\ell(\alpha-)-\phi(\alpha-)$,
adding and subtracting $\eta_\ar(\alpha-)$ to the right hand
side of (\ref{est:21}), we obtain
\begin{eqnarray*}
\eta'_\ell(\alpha) &=&
-c_0'+c_0+\eta_\ar'(\alpha-)-\psi(\alpha)+\psi(\alpha-)+\nu(\alpha)+\ell(\alpha)+\eta_\ell(\alpha-)
-\phi(\alpha-)-\eta_\ar(\alpha-).
\end{eqnarray*}
From (\ref{est:1}) we infer that $\eta'_\ar(\alpha-)  \leq
\eta_\ar(\alpha-)+[c_0'-c_0]^+$, and so
\begin{eqnarray}
\eta'_\ell(\alpha)&\leq &
-c_0'+c_0+[c_0'-c_0]^+-\psi(\alpha)+\psi(\alpha-)+\nu(\alpha)+\ell(\alpha)
\nonumber\\ & &+\eta_\ell(\alpha-)-\phi(\alpha-).
\label{est:2}
\end{eqnarray}
On the other hand, using the relations $\phi(\alpha)\geq \ell(\alpha)$ and
$\eta_r(\alpha)-\eta_r(\alpha-)\geq 0$, we have
\begin{eqnarray}\eta_\ell(\alpha)&=&
\eta_\ell(\alpha-)+\phi(\alpha)-\phi(\alpha-)+\psi(\alpha-)-\psi(\alpha)\label{est:3}\\
& &+\eta_\ar(\alpha)-\eta_\ar(\alpha-) \nonumber\\
&\geq &
\eta_\ell(\alpha-)+\ell(\alpha)-\phi(\alpha-)+\psi(\alpha-)-\psi(\alpha).\nonumber
\end{eqnarray} By combining (\ref{est:2}) and (\ref{est:3}), we
see that (\ref{est:4}) also holds, and the claim follows. A
similar argument shows that (\ref{est:4}) and (\ref{est:5}) are
also satisfied  when $\eta_\ar'(\alpha)-\eta_\ar'(\alpha-)>0$.

Next, note from the definition of $\alpha$ that there exists a
sequence of constants $\{s_n\}$ with $s_n\downarrow 0$ as
$n\rightarrow \infty$ such that one of the following statements
must be true:

(i)
$\eta_\ell'(\alpha+s_n)>\eta_\ell(\alpha+s_n)+\nu(\alpha+s_n)+[c_0-c_0']^+
\mbox{ for all }n \in \N$;

(ii) $\eta_\ar'(\alpha+s_n)> \eta_\ar(\alpha+s_n)+[c_0'-c_0]^+
\mbox{
  for all }n\in \N$.
\newline First, suppose Case (i) holds. Then,  taking the limit as
$n\rightarrow \infty$, by the right-continuity of $\eta_\ell$, $\eta'_\ell$
and $\nu$, we have $\eta'_\ell (\alpha) \geq \eta_{\ell} (\alpha) + \nu
(\alpha) + [c_0 - c_0']^+$.  Together with  (\ref{est:4}),  this implies
\begin{eqnarray}
\eta_\ell'(\alpha)=\eta_\ell(\alpha)+\nu(\alpha)+[c_0-c_0']^+.
\label{est:6}
\end{eqnarray}
Since $\eta_\ell$ and $\nu$ are non-decreasing, we have from
Case (i) and (\ref{est:6}) that
$\eta_\ell'(\alpha+s_n)>\eta_\ell'(\alpha)$ for each $n\in \N$.
By the complementarity condition in (\ref{comp-sp2}), this
implies $\phi'(\alpha)=\ell(\alpha)$. Together with
(\ref{est:5}), (\ref{est:6}) and the relation $\psi = \psi' +
\nu$, this implies
\begin{eqnarray*}
 \phi (\alpha) - \ell (\alpha)  =  \phi (\alpha) - \phi' (\alpha) & = &
c_0 - c'_0 + \nu (\alpha) + \eta_{\ell} (\alpha) - \eta'_{\ell} (\alpha) -
\eta_{\ar} (\alpha) + \eta'_{\ar} (\alpha) \\
& \leq & c_0 - c'_0 - [c_0 - c_0']^+ + [c_0' - c_0]^+ = 0.
\end{eqnarray*}
Since $\phi(\cdot) \in \gt$, this implies $\phi(\alpha) =
\ell(\alpha)$.

Consider the shift operator $T_{\alpha}:\dspace \mapsto
\dspace$ defined by $T_\alpha f (s) = f(\alpha + s) -
f(\alpha)$ for $s \in [0,\infty)$. By uniqueness of solutions
to the SP, it is easy to see that $(\phi(\alpha + \cdot),
T_{\alpha} \eta)$ solve the SP for $\phi(\alpha) + T_{\alpha}
\psi$ with the associated pair of constraining processes
  $(T_{\alpha} \eta_{\ell}, T_{\alpha} \eta_{\ar})$, and likewise for
$(\phi'(\alpha + \cdot), T_{\alpha} \eta')$. Now, by  the right
continuity of $\phi$, $\phi'$ and $\ar$ and the fact that
$\phi(\alpha) =\phi'(\alpha) = \ell(\alpha) < \ar(\alpha)$,
there exists $\varepsilon > 0$ such that for each $s\in
[0,\varepsilon]$, $\phi(\alpha+s)<r(\alpha+s)$ and
$\phi'(\alpha+s)<r(\alpha+s)$. The complementarity condition
(\ref{comp-sp2}) implies that $T_{\alpha} \eta = T_{\alpha}
\eta_{\ell}$ and $T_{\alpha} \eta' = T_{\alpha} \eta_{\ell}'$
on the interval $[0,\ve]$. An application of property (ii) of
Proposition \ref{lem:Mono}, with $c_0 = c_0' = \ell (\alpha)$
and $T_\alpha \psi$, $T_{\alpha} \psi'$, $\ell(\alpha+ \cdot)$,
$\ar(\alpha+\cdot)$ and $T_{\alpha} \nu$ in place of $\psi,
\psi', \ell, \ar$ and $\nu$, shows that for each $s\in
[0,\varepsilon]$, $T_{\alpha} \eta_{\ell}'(s) \leq T_{\alpha}
\eta_{\ell} (s) +T_{\alpha} \nu(s)$. 
Together with (\ref{est:6}), this shows that for each $s\in
[0,\varepsilon]$,
\[\eta_\ell'(\alpha+s)\leq
\eta_\ell(\alpha+s)+\nu(\alpha+s)+[c_0-c_0']^+, \]
which contradicts Case (i).
Hence Case (ii) should hold. In this case,  a similar argument can be used to
 show that $\phi'(\alpha)= \phi' (\alpha) = \ar(\alpha)$,
and  arguments analogous to those used above can then be applied to arrive at
a contradiction to Case (ii). Thus $\alpha=\infty$ or, in other
words, the second inequality in property 1 and the first
inequality in property 2 of the proposition hold. The first
inequality in property 1 and the second inequality in property 2
of the proposition can be proved in a similar way with $\beta$
instead of $\alpha$, where
\[\beta=\inf\{t> 0: \eta_\ar'(t)+\nu(t)+[c_0-c_0']^+<\eta_\ar(t) \mbox{  or }
\eta_\ell'(t)+[c_0'-c_0]^+<\eta_\ell(t)\}. \]
\end{proof}

\section{Variation of the Local Time of RBM}
\label{sec-rbm}

Throughout this section, let $B$ be a one-dimensional standard
BM starting from $0$ and defined on some filtered probability space
$(\Omega,\ \MF,\ \{\MF_t\},\ \P)$.  Also, let $\E$ denote expectation
with respect to $\P$.

\begin{defn}
Given $\ell \in \dspacemin$ and $\ar \in \dspacepl$ with $\ell \leq \ar$,
we define RBM $W$ on $[\ell (\cdot), \ar
(\cdot)]$ starting at $x \in \R$ by
\[ W = \bGamma_{\ell, \ar} (x + B). \]
\end{defn}

Due to uniqueness of solutions to the ESP, it is easy to see that
$\bGamma_{\ell, \ar}(\psi)(t)$ depends only on
$\{\ell(u), \ar(u)$, $\psi(u), u \in [0,t]\}$.
Thus, $W$ is adapted to the filtration generated by $B$.
 Moreover,  $W$ admits the unique
decomposition  $W(t)=x+B(t)+Y(t)$ for $t\geq 0$ such that for
each $\omega\in \Omega$,
$(W(\omega,\cdot),\ltime(\omega,\cdot))$ solves the ESP
 on $\gt$, as described in Definition \ref{def-smap2},
for $x +B(\omega, \cdot)$. We will refer to $Y$ as the local time of $W$ on the
(time-dependent) boundary of $\gt$.
{F}rom Corollary \ref{cor-separatelr}, it immediately follows
that $\ltime$  a.s.\ has finite  variation on every time
interval $[t_1,t_2]$ such that $\inf_{t\in[t_1,t_2]}(\ar(t) -
\ell(t)) > 0$.

{F}or the rest of this section, fix $\ell \in\dspacemin$, $\ar \in\dspacepl$ such that
$\ell \leq \ar$, define
 \bes
 \tau \doteq
 \inf\{t > 0:\ \ar(t)=\ell(t) \hbox {  or  }
 \ar(t-)=\ell(t-)\},
 \ees
and assume that $\tau \in (0,\infty)$. In Sections
\ref{subs-lower} and \ref{subs-upper} we identify some
necessary and some sufficient conditions for $\ltime$ to have
$\P$-a.s.\ finite variation on $[0,\tau]$. Recall that the
variation of a function $f$ on $[t_1,t_2]$ is denoted by
$\MV_{[t_1,t_2]}(f)$.   We apply these results in Section
\ref{subs-2drbm} to analyze the local time of a class of
two-dimensional RBMs in a fixed domain.

\subsection{A Lower Bound}
\label{subs-lower}

We show that the local time of RBM on $[0,\tau]$ has
infinite variation for some $\ell$ and $\ar$ by comparing the
space-time domain $\{(t,x): \ell (t) < x < \ar(t)\}$ to a
``comb domain.''

Let $K'$ denote a subset of $ \Z$, for example, $K'$ may be the
sequence of all negative integers, or all positive integers. We
denote by $K$ the subset of $K'$ consisting of all elements of
$K'$ except the largest element of $K'$, assuming one exists.

\begin{theorem}
\label{th-comb} Suppose that there exists a set $K'$ and a
sequence $\{s_k\}_{k\in K'}$ that is strictly increasing, takes
values in $[0,\tau]$ and, for some constant $c_1
\in(-\infty,\infty)$ and all $k\in K$, satisfies
 \be
 \label{eq:min}
 \frac{\min(\ar(s_{k+1}) - \ell(s_k),
 -\ell(s_{k+1}) + \ar(s_k))}{(s_{k+1} - s_k)^{1/2}}
 \leq c_1.
 \ee
If
 \be
 \label{eq:roots}
 \sum_{k\in K} (s_{k+1} - s_k)^{1/2} =\infty
 \ee
then $\MV_{[0,\tau]} \ltime =\infty$, a.s.
\end{theorem}

\begin{remark}
{\em
 The constant $c_1$ in the statement of Theorem
\ref{th-comb} does not have to be positive. Intuitively
speaking, the smaller $c_1$, the more the variation
accumulated by $\ltime$. Examples of domains that satisfy the
assumptions of the theorem are provided below the proof.
 }
\end{remark}

\begin{proof}[Proof of Theorem \ref{th-comb}]
Let $\Delta B_k = B(s_{k+1})- B(s_k) $, $c_2 = 1 \lor c_1$. We
define an event $A_k$ by
 \bes
 A_k =
 \{\Delta B_k \in( 2 c_2 (s_{k+1} - s_k)^{1/2},
 3c_2 (s_{k+1} - s_k)^{1/2})\}
 \ees
if $\ar(s_{k+1}) - \ell(s_k) \leq -\ell(s_{k+1}) + \ar(s_k)$,
and we let
 \bes
 A_k =
 \{\Delta B_k \in( -3 c_2 (s_{k+1} - s_k)^{1/2},
 -2 c_2 (s_{k+1} - s_k)^{1/2})\}
 \ees
if $\ar(s_{k+1}) - \ell(s_k) > -\ell(s_{k+1}) + \ar(s_k)$. By
Brownian scaling, there exists $p_1>0$ such that $\P(A_k) > p_1$
for all $k\in K$. This implies that
 \bes
 \E \left[|\Delta B_k| \ind_{A_k}\right] \geq p_1 2 c_2 (s_{k+1} - s_k)^{1/2},
 \ees
and so, in view of (\ref{eq:roots}), we have
 \bes
 \sum_{k \in K} \E\left[ |\Delta B_k| \ind_{A_k}\right] = \infty.
 \ees
We also have $|\Delta B_k| \ind_{A_k} \leq 3 c_2 (s_{k+1} -
s_k)^{1/2} \leq 3 \tau^{1/2}$, a.s., for every $k$. The random
variables $|\Delta B_k| \ind_{A_k}$ are independent. Hence, by
the ``three series theorem'' (\cite{durrett}, Ch.~1, (7.4)), we
have a.s.,
 \be
 \label{eq:sum}
 \sum_{k \in K}  |\Delta B_k| \ind_{A_k} = \infty.
 \ee

Suppose that the event $A_k$ holds and consider the case when
$\ar(s_{k+1}) - \ell(s_k) \leq -\ell(s_{k+1}) + \ar(s_k)$. Then
$W(s_k) \geq \ell(s_k) $ and $W(s_{k+1}) \leq \ar(s_{k+1})$.   Together
with (\ref{eq:min}) and the case assumption, this implies
$W(s_{k+1}) - W(s_k) \leq c_1(s_{k+1} - s_k)^{1/2}$. Since
$B(s_{k+1}) - B(s_k) \geq 2 c_2(s_{k+1} - s_k)^{1/2}$, we must
have $\ltime(s_{k+1}) - \ltime(s_k) \leq - c_2(s_{k+1} -
s_k)^{1/2}$. It follows that $\MV_{[s_k, s_{k+1}]} \ltime \geq
c_2(s_{k+1} - s_k)^{1/2} \geq (1/3) |\Delta B_k|$. A completely
analogous argument shows that the same bound holds in the case
when $\ar(s_{k+1}) - \ell(s_k) \geq -\ell(s_{k+1}) + \ar(s_k)$.
This estimate and (\ref{eq:sum}) imply that, a.s.,
 \bes
 \MV_{[0, \tau]} \ltime \geq
 \sum_{k \in K} \MV_{[s_k, s_{k+1}]} \ltime \geq
 \sum_{k \in K} (1/3) |\Delta B_k| \ind_{A_k} = \infty.
 \ees
\end{proof}

\begin{ex}
 {\em
It is straightforward to check that if $\ell(1/(2k)) \geq 0$
and $\ar(1/(2k+1)) \leq 0$ for $k \geq k_0$ then the
assumptions (\ref{eq:min}) and (\ref{eq:roots}) are satisfied
with $s_k = 1/k$, $k \geq 2k_0$, and the variation of the local
time is infinite on $[0,\tau]$, a.s.
 }
\end{ex}

\begin{ex}
\label{ex:1}
 {\em
Consider $\ell$ and $\ar$ such that $\ell(0) < \ar(0)$ and
$f(t) = \ar(t) - \ell(t)$ is a non-increasing function. Let the
sequence $\{s_k\}$ be defined in the following way. We let $s_0
= 0$, and for $k\geq 1$, we let $s_{k+1} = s_k + f^2(s_k)$.
Then $(s_{k+1} - s_k)^{1/2} = f(s_k)$ and (\ref{eq:min}) is
satisfied with $c_1 =1$, by construction. It follows that the
variation of the local time is infinite on $[0,\tau]$ a.s.,
provided $\sum_{k\geq 0} f(s_k) = \infty$.

Note that the number of intervals $[s_k, s_{k+1}]$ inside
$[\tau-2^{-j},\tau-2^{-j-1}]$ is bounded below by
$2^{-j-2}/f^2(\tau-2^{-j})$. The contribution of each one of
these intervals to the sum $\sum_k(s_{k+1} - s_k)^{1/2}$ is
bounded below by $f(\tau-2^{-j-1})$. Hence, the contribution
from all these intervals is bounded below by
$2^{-j-2}f(\tau-2^{-j-1})/f^2(\tau-2^{-j})$. It follows that if
for some $j_0$,
 \be
 \label{eq:ff2}
 \sum_{j> j_0}
 2^{-j-2}f(\tau-2^{-j-1})/f^2(\tau-2^{-j})
 =\infty
 \ee
then the variation of the local time is infinite on $[0,\tau]$
a.s.

Consider the case when $f(\tau-t) = t^\alpha$ for some
$\alpha>0$. If $\alpha \geq 1$ then (\ref{eq:ff2}) is true and
the variation of the local time is infinite on $[0,\tau]$ a.s.
 }
\end{ex}

\begin{ex}
\label{ex:2}
 {\em
This is a modification of the previous example. Suppose that
$\ell$ and $\ar$ are such that $\ell(0) = \ar(0)$ and $f(t) =
\ar(t) - \ell(t)$ is a non-decreasing function on some interval
$[0,\tau_1]$, with $\tau_1 \in(0,\tau)$. Assume that for some
$0<c_3,c_4<\infty$, we have $c_3<f(t)/f(2t)<c_4$ for all
$t\in(0,\tau_1/2)$. Let $\{s_k\}_{k\in K'}$ be the usual
ordering of all points of the form $2^{-j} + mf^2(2^{-j})$, for
$m=0,\dots, [2^{-j}/f^2(2^{-j})]$, and $j> j_1$, where $j_1$ is
chosen so that $s_k < \tau_1$ for all $k\in K'$. Then it is
easy to see that (\ref{eq:min}) is satisfied. Similarly, it is
routine to verify that the condition (\ref{eq:roots}) is
satisfied if $\sum_{j>j_1} 2^{-j} f(2^{-j}) =\infty$. Hence, if
$f(t) = t^\alpha$ for some $\alpha \geq 1$ and $t\in[0,\tau_1]$
then the variation of the local time is infinite on $[0,\tau]$
a.s.
 }
\end{ex}

\subsection{An Upper Bound}
\label{subs-upper}

Our  upper bound will be based on the comparison of the
space-time domain $\dot D\doteq\{(t,x): \ell (t) < x < \ar(t)\}$
with a family of ``parabolic boxes.''

\begin{theorem}
\label{th-boxes} Suppose that there exists a sequence
$\{s_k\}_{k\in \Z}$ that is strictly increasing and is such
that $\lim_{k\to -\infty} s_k = 0$, and $\lim_{k\to \infty} s_k
= \tau$. Suppose that there exist a constant $c_1 < \infty$ and
sequences $\{a_k\}_{k\in \Z}$ and $\{b_k\}_{k\in \Z}$, such
that $\{(t,x): s_k < t < s_{k+1}, a_k < x < b_k\} \subset \dot
D$, and
 \be
 \label{eq:boxx}
 (1/c_1) (s_{k+1} - s_k)^{1/2} < b_k - a_k < c_1 (s_{k+1} - s_k)^{1/2}
 \ee
for every $k\in \Z$. Further, given
$m_k \doteq (a_k+b_k)/2$,
\[ d_k \doteq |\ar(s_k) - m_k| + |\ell(s_k) - m_k| \]
and
\[ d_k' \doteq |\ar(s_{k+1}) - m_k| \vee |\ar(s_{k+1}-)-m_k|
+ |\ell(s_{k+1}) - m_k| \vee |
\ell (s_{k+1}-) - m_k|,
\]
suppose that
\be
\label{eq:roots1}
 \sum_{k\in \Z} (s_{k+1} - s_k)^{1/2} < \infty,
\ee
\be
\label{eq:ra1}
 \sum_{k\in \Z} d_k < \infty  \qquad \qquad \mbox{ and } \qquad \qquad \sum_{k \in \Z} d_k' < \infty.
\ee
Then $\E\left[\MV_{[0,\tau]} \ltime\right] <\infty$.
\end{theorem}
\begin{proof}
Let $\diff_k \doteq b_k-a_k$.
{F}or some fixed $k\in \Z$, we will estimate the expected amount
of local time generated on an interval $[s_k, s_{k+1}]$.
{F}irst, choose
$\eps_0 \in(0,s_{k+1}-s_k)$
such that for all $\ve \in (0,\ve_0)$,
\begin{equation}
\label{eps-small}
 \left|\ell (s_{k+1} - \ve)- m_k \right| + \left|\ar (s_{k+1} -
 \ve) - m_k \right| \leq 2 d_k'.
\end{equation}
 Let
$T_1 \doteq s_k$, and for $\ve \in (0,\ve_0)$, define
 \beas
 U_1 &\doteq& \inf\{t\geq T_1: W(t) = m_k\} \land (s_{k+1}-\eps), \\
 T_j &\doteq& \inf\{t\geq U_{j-1}: |W(t) - m_k| \geq \diff_k/4 \}
 \land (s_{k+1}-\eps), \quad j\geq 2,\\
 U_j &\doteq& \inf\{t\geq T_j: W(t) = m_k\}
 \land (s_{k+1}-\eps) , \quad j\geq 2.
 \eeas
{F}or each $j \geq 2$,   $W$ is away from the
upper and lower boundaries on $[U_{j-1},T_j]$, and so   we have
 \be \MV_{[U_{j-1}, T_j]} \ltime =0 \quad \quad \mbox{ for all } j \geq
 2. \label{est:10} \ee

We now consider intervals of the form $[T_j,U_j]$, $j \in \N$.
The elementary relation $Y=W-B$ yields the  bound \be
 \left|\ltime(U_j) - \ltime(T_j)\right|
 \leq  \left|W(U_j) - W (T_j)\right| + \sup_{t\in[T_j, U_j]} |B(t) -
 B(T_j)| \quad
\label{est:11} \ee for every $j \in \N$. Standard estimates
show that there exists $c_2 < \infty$ such that
 \begin{eqnarray}
\nonumber
 \E\left[\sup_{t\in[T_j, U_j]} |B(t) - B(T_j)|\right]
 & \leq &  \E\left[\sup_{t\in[s_k, s_{k+1}]} |B(t) - B(s_k)|\right] \\
&  \leq & c_2 (s_{k+1} - s_k)^{1/2}.
\label{est:8}
 \end{eqnarray}
For every $j \geq 1$ such that $U_j < s_{k+1} - \ve$, the
right-continuity of $W$ ensures that $W(U_j) = m_k$. Likewise,
for every $j \geq 2$ such that $T_j < s_{k+1} - \ve$, we have
$|W(T_j) - m_k| = \Delta/4$ because, as is easy to see, $W$ is
continuous at $T_j$. Indeed, the latter assertion follows
because $B$ is continuous, $W(t) \in [m_k - \Delta_k/4, m_k +
\Delta_k/4] \subset (\ell(t), \ar(t))$ for every $t \in  [T_j,
U_j]$ and, by equations (\ref{implication-baretazero}) and
(\ref{jump}), at any jump time $t$ of $W$, either $\ell (t-) =
W(t-), \ell (t) = W(t)$ or $\ar(t-) = W(t-), \ar(t) = W(t)$.
The last two statements, when combined with the triangle
inequality, the fact that $W(s) \in [\ell(s),\ar(s)]$ for every
$s$, and the relation (\ref{eps-small}),  show that
\begin{eqnarray}
\nonumber
|W(U_1) - W(T_1)| & \leq & \ind_{\{U_1 = s_{k+1} - \ve\}}
|W(U_1) - m_k| + |m_k - W(T_1)| \\
& \leq & 2d_k' + d_k
\end{eqnarray}
and, for $j \geq 2$,
\begin{eqnarray}
\nonumber
|W(U_j) - W(T_j)| & = & \ind_{\{T_j < s_{k+1} - \ve\}}|W(U_j) - W(T_j)| \\
\nonumber
& \leq & \ind_{\{T_j < s_{k+1} - \ve\}}[\ind_{\{U_j = s_{k+1} - \ve\}}
|W(U_j) - m_k| + |m_k - W(T_j)|] \\
& \leq &  \ind_{\{T_j < s_{k+1} - \ve\}} [2d_k' + \Delta_k/4].
\end{eqnarray}
In turn,  together with (\ref{est:11}) and (\ref{est:8}), this implies that
\begin{equation}
\label{eq-ltime1}
\E[|Y(U_1) - Y(T_1)|] \leq   2 d_k' + d_k +  c_2 (s_{k+1} - s_k)^{1/2}
\ee
and, for $j \geq 2$,
\begin{eqnarray}
\nonumber
\E[|Y(U_j) - Y(T_j)|] & \leq & \left( 2 d_k' + \Delta_k/4 +  c_2 (s_{k+1} -
s_k)^{1/2}\right) \E[\ind_{\{T_j < s_{k+1} - \ve\}}] \\
\label{eq-ltimej}
& \leq & \left( 2 d_k' +  c_3 (s_{k+1} -
s_k)^{1/2}\right) \E[\ind_{\{T_j < s_{k+1} - \ve\}}],
\end{eqnarray}
where the last inequality holds with $c_3 = c_2 + c_1/4$ due to
(\ref{eq:boxx}).

Now, if $T_j < s_{k+1} -\eps$ then the process $B$ must have
had an oscillation of size $\diff_k/4$ or larger inside the
interval $[U_{j-1}, T_j]$.  However, $\diff_k/4 \geq (s_{k+1} -
s_k)^{1/2}/4 c_1$ by  the inequality (\ref{eq:boxx}), and so it
can be deduced from the Kolmogorov-\v{C}entsov theorem that the
expected number of oscillations of $B$ of size $ (s_{k+1} -
s_k)^{1/2}/4 c_1$ on the time interval $[s_k, s_{k+1}]$ is
bounded by a constant $c_4<\infty$. In other words,
 \be
\ds \sum\limits_{j\geq 2}\E\left[ \ind_{\{T_j < s_{k+1} -\eps\}}\right]
 =
\E\left[\sum_{j\geq 2}  \ind_{\{T_j < s_{k+1} -\eps\}}\right]
 \leq c_4. \label{est:16}
 \ee
 Summing  (\ref{eq-ltimej}) over $j \geq 2$, adding (\ref{eq-ltime1}) and
using (\ref{est:16}), we obtain
\be
 \ds \sum\limits_{j\geq 1} \E\left[ \left|\ltime (U_j) - \ltime (T_j)
 \right|\right]
 \leq d_k + c_5 d_k'
+ c_5 (s_{k+1}-s_k)^{1/2},
\ee
where  $c_5= (c_3 c_4 + c_2) \vee (2 + 2 c_4)$.

Since
$W$ touches at most one boundary on each interval $[T_j, U_j]$, $j \in \N$,
 $Y$ is monotone on each such interval.  Thus
\[
\E \left[\MV_{[s_k, s_{k+1} - \eps]} \ltime\right]  \leq
\sum\limits_{j\geq 1} \E\left[ \left|\ltime (U_j) - \ltime (T_j) \right|
  \right]  \leq  d_k +  c_5 d_k'
+ c_5 (s_{k+1}-s_k)^{1/2}.
\]
By taking the limit as $\varepsilon\rightarrow 0$ and using the fact that
variation is monotone, we conclude that
 \be
\label{bd-var-}
 \E \left[\MV_{[s_k, s_{k+1})} \ltime\right]
 \leq d_k +  c_5 d_k'
+ c_5 (s_{k+1}-s_k)^{1/2}.
 \ee
The process $\ltime$ may have a jump at time $s_{k+1}$, whose
size can be bounded,  using the relation $Y = W-B$, the
triangle inequality, the continuity of the paths of $B$ and the
fact that $W(s) \in [\ell(s), \ar(s)]$ for all $s$, as follows:
 \beas
 |Y(s_{k+1}) - Y(s_{k+1}-)| & \leq &   |W(s_{k+1}) - W(s_{k+1}-)| \\
&  \leq &   |W(s_{k+1}) - m_k|  +   |W(s_{k+1}-) - m_k| \leq  2d_k'.
 \eeas
Together with (\ref{bd-var-}), this implies that
 \[
 \E \left[\MV_{[s_k, s_{k+1}]} \ltime\right]
\leq d_k + (c_5 + 2) d_k' + c_5 (s_{k+1}-s_k)^{1/2}.
\]
Summing over $k$ and using (\ref{eq:roots1}) and (\ref{eq:ra1}), we
obtain
 \bes
 \E \left[\MV_{[0,\tau]} \ltime\right]
 \leq \sum_{k\in \Z}
\left( d_k + (c_5 + 2) d_k' + c_5 (s_{k+1}-s_k)^{1/2} \right)
 < \infty.
 \ees
\end{proof}

\begin{ex}
 {\em
Our first example is elementary. Let $-\ell(t) = r(t) =
t^\alpha$ for $t\in [0, \tau/4]$ and $-\ell(t) = r(t) =
(\tau-t)^\alpha$ for $t\in [3\tau/4, \tau]$, where $\alpha >0$.
We assume that $\ell$ and $\ar$ are continuous on $[0,\tau]$
and $\ar(t) > \ell(t)$ for $t\in(0,\tau)$. Let $f(t) = \ar(t)
-\ell(t)$. Let $\{s_k\}_{k\in \Z}$ be the usual ordering of all
points belonging to two families: (i) all points of the form
$2^{-j} + mf^2(2^{-j})$, for $m=0,\dots, [2^{-j}/f^2(2^{-j})]$,
and $j> j_1$, where $j_1$ is the smallest integer such that
$2^{-j_1} < \tau/4$, and (ii) all points of the form $\tau -
2^{-j} - mf^2(\tau-2^{-j})$, for $m=0,\dots,
[2^{-j}/f^2(\tau-2^{-j})]$, and $j> j_1$. We let $a_k$ be the
smallest real number, and we let $b_k$ be the largest real
number such that $\{(t,x): s_k < t < s_{k+1}, a_k < x < b_k\}
\subset \dot D$. It is easy to verify that (\ref{eq:boxx})
holds.

Recall the notation from the proof of Theorem \ref{th-boxes}.
Consider an interval $[s_k, s_{k+1}] \subset [2^{-j},
2^{-j+2}]$. Then $d_k \lor d'_k \leq c_1 f(2^{-j})$ and
$(s_{k+1} - s_k)^{1/2} \leq f(2^{-j+1})$. The sum over all $k$
in the indicated range gives us
 \bes
 \sum (d_k \lor d'_k) + (s_{k+1} - s_k)^{1/2}
 \leq c_2 2^{-j}f(2^{-j+1})/f^2(2^{-j}).
 \ees
Summing over $j> j_1$ yields a finite number provided $\alpha
<1$. A similar analysis applies in the interval $[3\tau/4,
\tau]$, so the assumptions (\ref{eq:roots1}) and
(\ref{eq:ra1}) of
Theorem \ref{th-boxes} are satisfied if $\alpha < 1$. We
conclude that if $\alpha < 1$ then $\E\left[ \MV_{[0,\tau]} \ltime\right] <
\infty$.
 }
\end{ex}

\begin{ex}
\label{ex:3}
 {\em
We present a stronger version of the last example, in which we relax
the assumption of symmetry between $\ell$ and $\ar$, but impose a little
more regularity of the paths $\ell$ and $\ar$.  Let $f(t) =
\ar(t) -\ell(t)$, $f(t) = t^\alpha$ for $t\in [0, \tau/4]$ and
$f(t) = (\tau-t)^\alpha$ for $t\in [3\tau/4, \tau]$, where
$\alpha >0$. We assume that $\ar(t) > \ell(t)$ for $t\in(0,\tau)$. The
crucial assumption in this example is that both $\ell$ and
$\ar$ are H\"older continuous with some exponent $\beta > 1/2$,
i.e., for some $c_1<\infty$ and all $t_1,t_2\in [0,\tau]$, we
have $|\ell(t_1) - \ell(t_2)| \leq c_1 |t_1 - t_2|^\beta$, and
a similar formula holds for $\ar$.

We proceed as in the previous example. Let $\{s_k\}_{k\in \Z}$
be the usual ordering of all points belonging to two families:
(i) all points of the form $2^{-j} + mf^2(2^{-j})$, for
$m=0,\dots, [2^{-j}/f^2(2^{-j})]$, and $j> j_1$, where $j_1$ is
the smallest integer such that $2^{-j_1} < \tau/4$, and (ii)
all points of the form $\tau - 2^{-j} - mf^2(\tau-2^{-j})$, for
$m=0,\dots, [2^{-j}/f^2(\tau-2^{-j})]$, and $j> j_1$. We let
$a_k$ be the smallest real number, and we let $b_k$ be the
largest real number such that $\{(t,x): s_k < t < s_{k+1}, a_k
< x < b_k\} \subset \dot D$. One can verify, as in the previous
example, that (\ref{eq:boxx}) holds; this is more involved but
sufficiently straightforward that we leave it to the reader.

Our assumption that $\ell$ and $\ar$ are H\"older continuous
with an exponent greater than $1/2$ can be used to show that
$d_k \lor d'_k \leq c_1 f(2^{-j})$ for $k$ and $j$ such that
$[s_k, s_{k+1}] \subset [2^{-j}, 2^{-j+2}]$. We also have
$(s_{k+1} - s_k)^{1/2} \leq f(2^{-j+1})$.

The rest of the analysis proceeds as in the previous example.
The sum over all $k$ in the indicated range gives us
 \bes
 \sum_k \left[ (d_k \lor d'_k) + (s_{k+1} - s_k)^{1/2} \right]
 \leq c_2 2^{-j}f(2^{-j+1})/f^2(2^{-j}).
 \ees
Summing over $j> j_1$ yields a finite number provided $\alpha
<1$. A similar analysis applies in the interval $[3\tau/4,
\tau]$, and so the assumptions (\ref{eq:roots1}) and
(\ref{eq:ra1}) of
Theorem \ref{th-boxes} are satisfied if $\alpha < 1$. We
conclude that if $\alpha < 1$ then $\E\left[ \MV_{[0,\tau]} \ltime\right] <
\infty$.
 We see that within the family of functions $f(\cdot)$ that
decay towards the endpoints of $[0,\tau]$ as $t^\alpha$, our
results are sharp, by comparing the present example with
Examples \ref{ex:1} and \ref{ex:2}.
}
\end{ex}

\begin{remark} {\em
Note that the parameters $\alpha$ and $\beta$ in Example \ref{ex:3}
can be such that $1/2 < \beta < \alpha < 1$. Consider a
function $\ell$ that is H\"older continuous with exponent
$\beta$ but it is not H\"older continuous with exponent
$\beta+\eps$ on any interval $[0, s]$, for any $\eps>0$ and any
$s>0$. A typical trajectory of a fractional Brownian motion
with appropriate exponent provides an example of such function.
By making a linear transformation, we may assume that $\ell(0)=
\ell(\tau) = 0$. Let $f(t) = t^\alpha$ for $t\in [0, \tau/4]$
and $f(t) = (\tau-t)^\alpha$ for $t\in [3\tau/4, \tau]$, $f(t)$
is continuous on $[0,\tau]$ and $f(t)> 0$ for $t\in(0,\tau)$.
Let $\ar(t) = \ell(t) +f(t)$. Then $\ell$ and $\ar$ satisfy the
assumptions of the present example, so $\E \left[\MV_{[0,\tau]} \ltime\right]
< \infty$. Note that neither $\ell$ nor $\ar$ need be monotone,
and both functions can oscillate between positive and negative
values. Their local oscillations near $0$ may be comparable in
absolute value to $t^\beta$, a function much larger than $f(t)
= t^\alpha$.
 }
\end{remark}

\subsection{Analysis of a class of 2-dimensional RBMs}
\label{subs-2drbm}

We now apply the results obtained in the last two sections to analyze
a class of two-dimensional RBMs studied in \cite{burtoby} and \cite{wil} (see also \cite{ramrei1, ramrei2} to see
how RBMs in this class arise as diffusion approximations of a class of queueing networks),
which  provided one of the sources of motivation for the current work.

We begin by recalling the setup from \cite{burtoby} and rephrase some of the results from that paper
in our terminology. The domain $D \subset \R^2$ is described by two continuous real-valued
functions $L$ and $R$ defined on $[0,\infty)$ that satisfy $L(0) = R(0) = 0$
and $L(y) < R(y) $ for all $y>0$. Let $D$ be given by
 \bes
 D\equiv \{ (x,y) \in
 \R^2 : y\geq 0, L(y) \leq x \leq R(y) \}.
 \ees
Let $\partial D^1 = \{ (x,y)\in \partial D : x = L(y) \}$ and
$\partial D^2 = \{ (x,y)\in\partial D : x = R(y) \}$. The paper
\cite{burtoby} was concerned with the two-dimensional RBM  $Z= (Z^1, Z^2)$ in $D$, with the vectors of
reflection horizontal on $\partial D^1$ and $\partial D^2$,
and an additional vertical direction of reflection at $0$ that ensures the RBM stays within $D$ (see
\cite{burtoby} for details).
{F}rom the Skorokhod-type lemma proved in \cite{burtoby}, it follows that $Z^2$ is a one-dimensional
RBM on $[0,\infty)$ and $Z^1$ is the ESM applied to a standard $1$-dimensional
  BM $B$ in the domain with time-dependent boundaries $\ell(t) = L(Z^2(t))$ and
$\ar(t) = R(Z^2(t))$.  By Definition \ref{def-smap2}, $Z^1$ admits the
decomposition $Z^1 = B + Y$ where $Y$ is (pathwise) the local time
or the pushing term associated with the ESP on $\gt$.

We first study the total variation of the local time $Y$ on a single excursion of
the RBM $(Z^1, Z^2)$ from the origin.
Let $[\tau_1,\tau_2]$ be an excursion interval for $Z^2$, i.e.,
$\tau_1 < \tau_2$, $Z^2(\tau_1) = Z^2(\tau_2) = 0$, and $Z^2(t)
> 0$ for $t\in(\tau_1, \tau_2)$.
The  following result was established in Theorem 3 of \cite{burtoby} -- we provide
an alternative proof of this result.

\begin{prop}
\label{lem-burtoby} Suppose that there exist $\varepsilon >0$
and $\gamma>2$ such that $R(y) - L(y)\leq y^\gamma$ for $y
\in[0,\eps]$. Then
\[ \MV_{[\tau_1, \tau_2]} \ltime = \infty   \quad a.s. \]

On the other hand, suppose that there exist $\varepsilon >0$
and $\gamma < 2$ such that $R(y) - L(y)\geq y^\gamma$ for $y
\in[0,\eps]$, and $R$ and $L$ are Lipschitz. Then
\[ \MV_{[\tau_1, \tau_2]} \ltime < \infty   \quad a.s. \]
\end{prop}
\begin{proof}
We start with the first case.  Let $\gamma_1 <1/2$
be such that $\gamma \cdot \gamma_1 >1$. Path properties at
endpoints of an excursion of (1-d reflected) Brownian motion from 0 are well known
to be the same as those of the 3-dimensional Bessel process, see,
e.g., \cite{burbook}. Hence, it follows from Theorem 3.3 (i) of
\cite{shigawat} that
 \bes
 \limsup_{t\downarrow \tau_1}
 \frac{Z^2(t-\tau_1)}{ (t-\tau_1)^{\gamma_1}}
 =0.
 \ees
This implies that
 \bes
 \limsup_{t\downarrow \tau_1} \frac
 {R(Z^2(t-\tau_1))-L(Z^2(t-\tau_1))}
 { (t-\tau_1)^{\gamma_1\gamma}}
 =0.
 \ees
If we set $\tau = \tau_2 - \tau_1$, $\ell(t) =
L(Z^2(t-\tau_1))$, $\ar(t) = R(Z^2(t-\tau_1))$ and $f(t) = \ar(t) - \ell(t)$,
then we see that $f(t) \leq t^{\gamma_1\gamma}$
for $t$ sufficiently close to $0$, where
$\gamma_1\gamma >1$. It follows from Theorem \ref{th-comb} that
the variation of local time accumulated by $Z^1$ on the
interval $[\tau_1,\tau_2]$ is infinite a.s. (see also Example
\ref{ex:2}).

Next, suppose that $\gamma < 2$, there exists $\varepsilon >0$ such
that $R(y) - L(y)\geq y^\gamma$ for $y \in[0,\eps]$, and the
functions $L$ and $R$ are Lipschitz. Let $\gamma_2> 1/2$ be such
that $\gamma_2 \gamma<1$. We use Theorem 3.3 (ii) of
\cite{shigawat} to see that
 \bes
 \liminf_{t\downarrow \tau_1}
 \frac{Z^2(t-\tau_1)}{ (t-\tau_1)^{\gamma_2}}
 =\infty.
 \ees
It follows that
 \bes
 \liminf_{t\downarrow \tau_1} \frac
 {R(Z^2(t-\tau_1))-L(Z^2(t-\tau_1))}
 { (t-\tau_1)^{\gamma_2\gamma}}
 =\infty.
 \ees
Using the notation introduced above, $f(t) \geq
t^{\gamma_2\gamma}$ for $t$ sufficiently close to $0$, where
$\gamma_2\gamma <1$. We would like to apply Theorem
\ref{th-boxes}. We proceed with the construction of boxes as in
Example \ref{ex:3}. The only new subtle point in the argument
is the verification of (\ref{eq:ra1}).
In Example \ref{ex:3}, we used
the fact that $\ell$ and $\ar$ were H\"older continuous with
exponent $\beta > 1/2$. Now, we use the fact that for any $\beta_1 \in (0,1/2)$, Brownian
motion is H\"older continuous with exponent $\beta_1$, a.s.,
and that the same applies to the trajectories of the 3-dimensional
Bessel process (because Brownian excursions from 0 have the
same local path properties). Since $L$ and $R$ are assumed to
be Lipschitz, we conclude that $\ell$ and $\ar$ are H\"older
continuous with  some exponent $\beta_1 >
\gamma_2\gamma/2$. This suffices to prove that the inequalities in
(\ref{eq:ra1}) hold. A similar analysis applies
at the other endpoint of the excursion, i.e., close to
$\tau_2$. We conclude that  the variation of the local time
accumulated by $Z^1$ on the interval $[\tau_1,\tau_2]$ is
finite a.s.
\end{proof}

We  now consider a somewhat different, and perhaps more
natural, question of whether $Z$ is a semimartingale. A
surprising fact is that for any functions $R$ and $L$ such that
$R(0) = L(0)$, the amount of local time accumulated on the
boundary is infinite. The rate of growth of $R - L$ in a
neighborhood of 0 turns out to be irrelevant for this question.
Our argument is based exclusively on the scaling properties of
Brownian motion. The proof will use excursion theory; see
\cite{burbook} for a review of the relevant definitions and
facts.

\begin{prop}
Suppose $Z(0)=0$. Then for every $T>0$,
\[  \MV_{[0,T]} \ltime = \infty, \quad \text{a.s.} \]
Consequently, the process $Z$ starting from $0$ is not a
semimartingale.
\end{prop}
\begin{proof}
 Let $[s_k,t_k]$, $k\geq 1$, be the collection of all excursion intervals
of $Z^2$ from 0. In other words, we have $Z^2(s_k) = Z^2(t_k)
=0$ and $Z^2(t) > 0$ for $t\in (s_k, t_k)$. Let $\{\sigma(t), t \geq 0\}$ be
the local time of $Z^2$ at 0. Then the family $\{(\sigma ({s_k}),
\{Z^2(t), t \in [s_k,t_k]\})\}_{k\geq 1}$ is a Poisson point
process on the space $[0,\infty) \times {\cal U}$, where $\cal
U$ is the space of excursions. The intensity of the Poisson
point process is the product of Lebesgue measure and an
excursion law $H$.

Note that we have $Z^1 = B + Y$ and $B$ is independent of
$Z^2$. Let $u_k \in [s_k,t_k]$ be the time when $Z^2$ attains
its maximum on the interval $[s_k,t_k]$. Given $s_k$ and $t_k$,
the process $\{(B(t) - B(u_k),Z^2(t)-Z^2(u_k)), t\in [u_k,
t_k]\}$ is independent of all processes
$\{(B(t)-B(u_j),Z^2(t)-Z^2(u_j)), t\in [u_j, t_j]\}$, $j\ne k$.

Given $t_k - s_k = a$, the distribution of $a^{-1}(t_k-u_k)$ is
independent of $a$, by scaling. Given $t_k - s_k = a$, the
distribution of $a^{-1/2}(B(t_k) - B(u_k))$ is otherwise
independent of $t_k$ and $s_k$, and of $Z^1(u_k)$. By Brownian
scaling, there exists  $c_0>0$ such that
 \bes
 \P(B(t_k) - B(u_k) > c_0 a^{1/2} \mid
 t_k - s_k = a) \geq 1/4,
 \ees
and
 \bes
 \P(B(t_k) - B(u_k) < - c_0 a^{1/2} \mid
 t_k - s_k = a) \geq 1/4.
 \ees
Using independence from $Z^1(u_k)$,
 \bes
 \P(Z^1(u_k)+ B(t_k) - B(u_k) > c_0 a^{1/2} \mid
 t_k - s_k = a, Z^1(u_k) >0 ) \geq 1/4,
 \ees
and
 \bes
 \P(Z^1(u_k)+ B(t_k) - B(u_k) < - c_0 a^{1/2} \mid
 t_k - s_k = a, Z^1(u_k) <0 ) \geq 1/4.
 \ees
Combining the two cases,
 \bes
 \P(|Z^1(u_k)+B(t_k) - B(u_k)| > c_0 a^{1/2} \mid
 t_k - s_k = a) \geq 1/4.
 \ees
Note that $Y(u_k) - Y(t_k) = Z^1(u_k)+B(t_k) - B(u_k)$. Since
$Z(t_k)=0$,
 \bes
 \P(|Y(u_k) - Y(t_k)| > c_0 a^{1/2} \mid
 t_k - s_k = a) \geq 1/4.
 \ees
This implies that
 \be\label{exc1}
 \P(\MV_{[u_k,t_k]} Y > c_0 a^{1/2} \mid
 t_k - s_k = a) \geq 1/4.
 \ee

Recall that $H$ is the excursion law for excursions of $Z^2$
from 0 and let $\zeta$ be the lifetime of an excursion. Then
$H(\zeta \in da) = c_1 a^{-3/2}$ (see \cite{burbook}). It
follows from  excursion theory that the number of excursions
starting at a point $s_k \leq \sigma^{-1}(1)$ and such that
$t_k-s_k \in (2^{-j-1}, 2^{-j}]$ has the Poisson distribution
with the average $c_1 \int_{2^{-j-1}}^{2^{-j}} a^{-3/2} da =
c_2 2^{j/2}$. By (\ref{exc1}), the number of such excursions
with the property that $\MV_{[u_k,t_k]} Y > c_0 2^{-(j+1)/2}$
is minorized by the Poisson distribution with the average
$(1/4) c_2 2^{j/2}$. Hence $\sum_{s_k \leq \sigma^{-1}(1),
t_k-s_k \in (2^{-j-1}, 2^{-j}]} \MV_{[u_k,t_k]} Y$ is minorized
by a random variable which is the product of $c_0 2^{-(j+1)/2}$
and a Poisson random variable with the average $(1/4) c_2
2^{j/2}$. By excursion theory, the sums $\sum_{s_k \leq
\sigma^{-1}(1), t_k-s_k \in (2^{-j-1}, 2^{-j}]} \MV_{[u_k,t_k]}
Y$ are independent for different $j$. Now it is elementary to
check that, a.s.,
 \bes
  \MV_{[0,\sigma^{-1}(1)]} Y \geq
 \sum_{j\geq 1} \sum_{s_k \leq \sigma^{-1}(1), t_k-s_k \in (2^{-j-1},
 2^{-j}]} \MV_{[u_k,t_k]} Y = \infty.
 \ees
The same argument shows that $ \MV_{[0,\sigma^{-1}(t)]} Y =
\infty$, a.s., for every $t>0$. Thus $ \MV_{[0,s]} Y = \infty$,
a.s., for every $s>0$.

It is intuitively clear from the first part of the proof that
$Z^1$ is not a semimartingale. A subtle technical difficulty is
that it is not obvious that the Doob decomposition of $Z^1$ and
the Skorokhod representation have to be identical. We shall
show that this is indeed the case in the following paragraph.

Suppose that $Z^1$ is a semimartingale with the decomposition
$Z^1=\tilde B+\tilde Y$. Fix $k, m \in\N$, $k < m$, and define
a sequence of stopping times as follows.
Let $T_1^{k,m} =\inf\{t\geq 0: Z^2= 2^{-k}\}$, $S_1^{k,m}
 =\inf\{t\geq T_1^{k,m}: Z^2=2^{-m}\}$, and recursively
define, for each $n \geq 2$, $T_n^{k,m} =\inf\{t\geq
S_{n-1}^{k,m}: Z^2= 2^{-k}\}$ and $S_n^{k,m}=\inf\{t\geq
T_n^{k,m}: Z^2=2^{-m}\}$. On each interval
$[T_n^{k,m},S_n^{k,m}]$, $Z^1$ is  a semimartingale with
decomposition $\tilde{B} +\tilde{Y}$.
 On the other hand, since $\inf_{t\in
[T_n^{k.m},S_n^{k,m}]} R(Z^2(t))-L(Z^2(t))>0$,  by the last
assertion of Theorem \ref{th-sptv} and the uniqueness of Doob's
decomposition, it follows that a.s.\ $\tilde Y(t)-\tilde
Y(s)=Y(t)- Y(s)$ for all $s,t \in [T_n^{k,m}, S_n^{k, m}]$. In
particular,  this implies that, a.s., $\sum_{n=1}^\infty
\MV_{[T_n^{k,m},S_n^{k,m}]}Y=\sum_{n=1}^\infty
\MV_{[T_n^{k,m},S_n^{k,m}]}\tilde Y$. Letting $m \rightarrow
\infty$, and then $k \rightarrow \infty$, we have almost
surely, for each $r>0$,
\begin{eqnarray*}
\MV_{[0,\sigma^{-1}(r)]}\tilde Y  & \geq&  \lim_{k\rightarrow \infty}
\lim_{m\rightarrow \infty}\sum\limits_{n=1}^\infty \MV_{[T_n^{k,m}\wedge
\sigma^{-1}(r),S_n^{k,m} \wedge \sigma^{-1}(r)]}\tilde Y \\
&  = &  \lim_{k \rightarrow \infty}
\lim_{m \rightarrow \infty} \sum\limits_{n=1}^\infty \MV_{[T_n^{k,m} \wedge
\sigma^{-1}(r),S_n^{k,m} \wedge \sigma^{-1}(r)]} Y  \\
& \geq&   \sum_{j\geq 1} \sum_{s_k \leq \sigma^{-1}(r), t_k-s_k \in (2^{-j-1},
 2^{-j}]} \MV_{[u_k,t_k]} Y.
\end{eqnarray*}
 The first part of the proof now shows that
the last term equals infinity.
We conclude that $Z^1$, and
therefore $Z$, is not a   semimartingale.
\end{proof}

\bibliographystyle{amsplain}
\bibliography{refs4}

\end{document}